\documentclass{aims}
\usepackage{amsmath}
  \usepackage{paralist}
  \usepackage{graphics} 
  \usepackage{epsfig} 
\usepackage{graphicx}  
 \usepackage[colorlinks=true]{hyperref}
\hypersetup{urlcolor=blue, citecolor=red}

\def\currentvolume{X}

     \def\DOI{10.3934/xx.xx.xx.xx}

\newcommand{\sysmat}{A}

\newcommand{\im}{x}

\newcommand{\imorig}{\im_{\text{orig}}}

\newcommand{\imopt}{\im^\star}

\newcommand{\sino}{b}

\newcommand{\norm}[1]{\|#1\|}

\newcommand{\normone}[1]{\norm{#1}_1}

\newcommand{\normtwo}[1]{\norm{#1}_2}

\newcommand{\imsize}{N}
\newcommand{\imside}{\imsize_{\mathrm{side}}}

\newcommand{\sinosize}{M}

\newcommand{\pen}{J}

\def\nv{N_\text{v}}

\DeclareMathOperator*{\argmin}{argmin}

\newcommand{\ellone}{{\bf L1}}
\newcommand{\elloneineq}{\ellone$_\ineqpar$}

\newcommand{\spikes}{{\bf spikes}}

\newcommand{\ppower}{{\bf ${\boldsymbol p}$\,-power}}
\newcommand{\onepower}{{\bf 1-power}}
\newcommand{\twopower}{{\bf 2-power}}

\newcommand{\nrecov}{\nv^\ellone}
\newcommand{\nvsuf}{\nv^{\mathrm{suf}}}

\newcommand{\sparsity}{k}
\newcommand{\relsparsity}{\kappa}
\newcommand{\relsampling}{\mu}

\newcommand{\emdash}{\,---\,}

\newcommand{\wone}{0.32}
\newcommand{\wtwo}{0.35}

\newcommand{\ineqpar}{\epsilon}
\newcommand{\numthr}{\tau}

\title[Relation between sampling and sparsity in CT]
      {Empirical average-case relation between undersampling and sparsity in x-ray CT}

\author[J S J\o{}rgensen, E Y Sidky, P C Hansen and X Pan]{}

\subjclass{Primary: 90C90, 15A29, 44A12; Secondary: 94A08.}
 \keywords{Inverse Problems, Computed Tomography, Image Reconstruction, Compressed Sensing, Sparse Solutions.}

\begin{document}
\maketitle

\centerline{\scshape Jakob S. J\o{}rgensen}
\medskip
{\footnotesize
 \centerline{Department of Applied Mathematics and Computer Science, Technical University of Denmark}
   \centerline{Richard Petersens Plads, Building 324, 2800 Kgs. Lyngby, Denmark, \url{jakj@dtu.dk}}
} %

\medskip

\centerline{\scshape Emil Y. Sidky}
\medskip
{\footnotesize
 \centerline{Department of Radiology, University of Chicago}
   \centerline{5841 South Maryland Avenue, Chicago, IL 60637, USA, \url{sidky@uchicago.edu}}
}

\medskip

\centerline{\scshape Per Christian Hansen}
\medskip
{\footnotesize
 \centerline{Department of Applied Mathematics and Computer Science, Technical University of Denmark}
   \centerline{Richard Petersens Plads, Building 324, 2800 Kgs. Lyngby, Denmark, \url{pcha@dtu.dk}}
} %

\medskip

\centerline{\scshape Xiaochuan Pan}
\medskip
{\footnotesize
 \centerline{Department of Radiology, University of Chicago}
   \centerline{5841 South Maryland Avenue, Chicago, IL 60637, USA, \url{xpan@uchicago.edu}}
}

\bigskip

 \centerline{\today}

\begin{abstract}
In x-ray computed tomography (CT) it is generally acknowledged that reconstruction
methods exploiting image sparsity allow reconstruction from a significantly reduced
number of projections. 
The use of such reconstruction methods is motivated by recent progress in compressed sensing (CS).
However, the CS framework provides neither guarantees of accurate CT reconstruction,
nor any relation between sparsity and a sufficient number of measurements for recovery, i.e., perfect reconstruction from noise-free data.
We consider reconstruction through 1-norm minimization, as proposed in CS, from data obtained using a standard CT fan-beam sampling pattern.
In empirical simulation studies we establish quantitatively a relation between the image sparsity and the
sufficient number of measurements for recovery within image classes motivated by tomographic applications.
We show empirically that the specific relation depends on the image class and in many cases exhibits a sharp phase transition as seen in CS, i.e. same-sparsity image require the same number of projections for recovery. Finally we demonstrate that the relation holds independently of image size 
and is robust to small amounts of additive Gaussian noise.
\end{abstract}

\section{Introduction}

In x-ray computed tomography (CT) an image of an object is reconstructed from projections
obtained by measuring the attenuation of x-rays passed through the object.
Motivated by reducing the exposure to radiation, there is a
growing interest in low-dose CT, cf. \cite{Yu:2009} and references therein.
This is relevant many fields, for example, in medical imaging to reduce the risk of radiation-induced
cancer, in biomedical imaging where high doses can damage the
specimen under study, and in materials science and non-destructive testing to reduce scanning time.

Classical reconstruction methods are based on closed-form analytical or approximate inverses of the
continuous forward operator;
examples are the filtered back-projection method \cite{natterer1986mathematics} and
the Feldkamp-Davis-Kress method for cone-beam CT \cite{FDK84}.
Their main advantages are low memory demands and computational
efficiency, which make them the current methods of choice in
commercial CT scanners \cite{PanIP:09}. However, they are known to have limitations on reduced data.

Alternatively, an algebraic formulation can be used, in which both object and data domains are discretized, 
yielding a large sparse system of linear equations.
This approach can handle non-standard scanning geometries for which no analytical inverse is available.
Furthermore, data reduction arising from low-dose imaging
can sometimes be compensated for by exploiting prior information about the image, such as smoothness or as in our case sparsity, i.e., having a representation with few non-zero coefficients. 

Developments in compressed sensing (CS) \cite{candes2006robust, Donoho2006} show potential
for a reduction in data while maintaining or even improving reconstruction quality.
This is made possible by exploiting image sparsity; loosely speaking, if the image
is ``sparse enough'', 
it admits accurate reconstruction from undersampled data. 
We refer to such methods as \emph{sparsity-exploiting methods}.
Under specific conditions on the sampling procedure, e.g., incoherence, there exist \emph{guarantees} of perfect recovery

Different types of sparsity are relevant in CT\@. For blood vessels \cite{li2002accurate} the image itself can be considered sparse, and reconstruction based on minimizing the $1$-norm can be expected to work well.  
The human body consists of well-separated areas of relatively homogeneous tissue and many materials, such as metals, consist of non-overlapping uniform sub-components. In these cases the image gradient is approximately sparse, and reconstruction based on minimizing total variation (TV)  of the image \cite{rudin1992nonlinear} is often a good choice.
Empirical studies using simulated as well as clinical data, using standard highly structured, i.e., non-random, sampling patterns, have
demonstrated that sparsity-exploiting methods indeed allow for reconstruction
from fewer projections
\cite{Bian:10,han2011algorithm,ritschl2011improved,SidkyTV:06,sidky2008image}.

However, as will be explained in Section \ref{sec:guarantees}, CS guarantees of recovery from undersampled data fall far short of extending to the sampling done in CT.
Therefore, in spite of the positive empirical results, we still lack a fundamental understanding
of conditions\emdash{}especially in a quantitative sense\emdash{}under which such methods can be expected to perform well in CT\@: How sparse images can be reconstructed and from how few samples?

The present paper demonstrates empirically an average-case relation between image sparsity 
and the sufficient number of CT projections enabling image recovery.
Inspired by work of Donoho and Tanner~\cite{DonohoTanner:2009} we use computer simulations 
to empirically study \emph{recoverability} within well-defined classes of sparse images.
Specifically, we are interested in the average number of CT projections sufficient for
exact recovery of an image as function of the image sparsity.
To simplify our analysis we focus on images with sparsity in the image domain and reconstruction 
through minimization of the image 1-norm subject to a data equality constraint, as motivated by CS.
These studies 
set the stage for forthcoming studies of other regularizers,
such as TV, as well as other types of sparsity.
We believe that our findings shed light on the connection between sparsity and sufficient sampling in CT and that they suggest the existence of a yet unknown theoretical foundation for CS in CT. 

The paper is organized as follows.
Section~\ref{sec:sparsityexploiting} presents the reconstruction problem of interest, describes relevant previous work and results from CS and discusses the application to CT.
Section~\ref{sec:experimentaldesign} describes all aspects of the experimental design, including the CT imaging model, generation of sparse images, and how to robustly solve the optimization problem of reconstruction.
Section~\ref{sec:numericalresults} contains our results establishing quantitatively a relation between image sparsity and sufficient sampling for recovery. Section~\ref{sec:noise} presents results for the noisy case and is followed by a discussion in Section~\ref{sec:discussion}.

\section{Sparsity-exploiting reconstruction methods}
\label{sec:sparsityexploiting}

This section describes our choice of sparsity-exploiting reconstruction method and explains how existing recovery guarantees from CS do not prove useful in the setting of CT.

\subsection{Reconstruction based on 1-norm minimization}

We consider the discrete inverse problem of recovering a signal
$\imorig \in \mathbb{R}^\imsize$ from data 
$\sino \in \mathbb{R}^\sinosize$.
The imaging model, which is assumed to be linear and discrete-to-discrete \cite{Barrett:FIS},
relates the image and the data through a \emph{system matrix}
$\sysmat \in \mathbb{R}^{\sinosize \times \imsize}$,
\begin{equation}\label{eq:eqsystem}
 \sysmat \im = \sino,
\end{equation}
where the elements of
$\im \in \mathbb{R}^\imsize$ are pixel values stacked into a vector. 
We focus in the present work on the ``undersampled'' and consistent case where $M<N$ such that \eqref{eq:eqsystem} has infinitely many solutions. To determine a unique solution we consider the problem
\begin{align} \label{eq:L1}
\text{\ellone{}:}\qquad\qquad \min_\im  \ \normone{\im} \quad
 \text{s.t.}\quad \sysmat \, \im = \sino ,
\end{align}
which seeks to recover the most sparse solution, i.e., the one with fewest non-zero components.
More generally, one could consider other forms of regularization based on prior knowledge or assumptions to restrict the set of solutions. 
A common type of regularization takes the form:
\begin{equation} \label{eq:ineqgeneralform}
\min_\im\ \pen(\im)\quad  \text{s.t.} \quad \normtwo{\sysmat\,\im - \sino} \leq \ineqpar,  
\end{equation}
where $\pen(\im)$ is the regularizer. The regularization parameter $\ineqpar$ controls the level of regularization and in the limit $\ineqpar \rightarrow 0$ we obtain the equality-constrained problem. 

The inequality-constrained problem \eqref{eq:ineqgeneralform} is of more practical interest than \eqref{eq:L1} because
it allows for noisy and inconsistent measurements, but its solution
depends in a complex way on the noise and inconsistencies in the data, as well
as the choice of the parameter~$\ineqpar$.
Studies of the equality-constrained problem, 
on the other hand,
provide an basic understanding of a given regularizer's reconstruction potential,
independent of specific noise.
Furthermore, many of the initial CS recovery guarantees deal with the equality-constrained formulation, and as such the equality-constrained problem in CT is the most natural place to first attempt to establish recovery results in CT.
Therefore, we focus for the major part of the present work on
the equality-constrained problem, however we conclude with a brief 
study of robustness with respect to additive Gaussian white noise.

\subsection{Existing recovery guarantees do not apply to CT} \label{sec:guarantees}

We call a tuple ($\imorig$, $\sysmat$) a \emph{problem instance} and say that $\imorig$
is \emph{recoverable} if the \ellone{} solution with $\sino = \sysmat\,\imorig$  is identical to $\imorig$. A caveat is that \ellone{} does not necessarily have a unique solution. 
The solution set may consist of a single image or an entire hyperface or hyperedge on the \hbox{$1$-norm} ball in case this set is aligned with the affine space of feasible points $\{ \im ~\vert~ \sysmat\im = \sino\}$. For being recoverable we therefore require that $\imorig$ be the \emph{unique} \ellone{} solution.

CS establishes \emph{guarantees} of recovery of sparse signals from a small number of measurements. An important concept is the restricted isometry property (RIP):
A vector $\im \in \mathbb{R}^{\imsize}$ with $\sparsity$ non-zero elements is called 
 $\sparsity$-sparse. A matrix $\sysmat$ is satisfies the RIP of order $\sparsity$ if a constant $\delta_\sparsity \in (0,1)$ exists such that 
\begin{equation}
 (1-\delta_\sparsity) \|\im\|_2^2 \leq \|\sysmat \im\|_2^2 \leq (1+\delta_\sparsity) \|\im\|_2^2.
\end{equation}
for all $\sparsity$-sparse vectors $\im$.
If the RIP-constant $\delta_\sparsity$ is small enough, then recovery of sparse signals is possible; more precisely, if $\delta_{2\sparsity} < \sqrt{2} -1$ then the \ellone{} solution $\imopt$ for data $\sino = \sysmat \imorig$ recovers the original image $\imorig$ \cite{candes2006stable}.
There also exist RIP-based guarantees for TV \cite{NeedellWard:SIR:13}. Such results are sometimes referred to as \emph{uniform} or \emph{strong} recovery guarantees \cite{FoucartRauhut:2013} since recovery of \emph{all} vectors of a given sparsity is ensured.

Unfortunately, the RIP is impractical to use because computing the RIP-constant for a given matrix is NP-hard 
\cite{Tillmann2014}. RIP-constants have been established only for few, very restricted class of matrices, e.g., with Gaussian i.i.d. entries \cite{candes2006robust}.
In \cite{SidkyPC:10} the authors computed lower bounds close to $1$ on the RIP-constant for very low sparsity values for x-ray CT system matrices, thereby excluding the possibility of RIP-based guarantees for tomography for other than extremely sparse signals. 

Other CS guarantees rely on incoherence \cite{CandesRomberg2007}, the spark \cite{Donoho2003spark} or the null-space property \cite{Cohen2009} but also fail to give useful guarantees, either due to being extremely pessimistic or NP-hard to compute, see, e.g., \cite{Dossal:2010}.
In fact, it is possible to construct examples of very sparse vectors that cannot be recovered from CT measurements \cite{PustelnikEUSIPCO2012}, implying that we cannot hope for uniform recovery guarantees for CT.

Instead, recoverability can be studied in an average-case sense. This can be motivated by a desire to ensure recovery of ``typical'' images of a given sparsity, but not all possible pathological images. Such results are sometimes referred to as \emph{non-uniform} or \emph{weak} recovery guarantees  \cite{FoucartRauhut:2013}.
Donoho and Tanner \cite{DonohoTanner:2009} derived weak recovery guarantees, i.e., established a relation between image sparsity and the critical average-case sampling level for recovery for example for Gaussian matrices. Their so-called phase diagram demonstrated agreement of the theoretical sampling level for a given sparsity with empirical recovery experiments by revealing a sharp phase transition between recovered and non-recovered images.

Weak recovery guarantees have been established for \emph{discrete} tomography \cite{PetraSchnoerr2014} using very restrictive sampling patterns (essentially only rays perfectly aligned with pixels) and for \emph{discrete-valued} tomography \cite{Gouillart2013}. These results, however, do not apply to the general case of x-ray CT, see \cite{HermanKuba:1999} for background on discrete tomography.
To our knowledge recovery guarantees for x-ray CT remain an open question.

\subsection{Our contribution}

In the present paper we empirically establish a quantitative relation between the number
of measurements and the image sparsity sufficient for average-case recovery. We generate random images from specific classes of images motivated by tomographic applications and determine the average critical sampling level for recovery by \ellone{} as function of image sparsity. 
Our empirical study is inspired by the Donoho-Tanner (DT) phase diagram, however we use a slightly modified diagram in which the quantities of our interest, namely sparsity and sampling, can be read directly off the axes.

With this approach we provide extensive empirical evidence of an average-case relation between sampling and sparsity occurring across different image classes, image sizes as well as showing robustness to noise.

\section{Design of numerical experiments}
\label{sec:experimentaldesign}
In this section we describe our overall empirical study design after presenting the chosen imaging model, generation of sparse test images and numerical optimization for the reconstruction problem.

\subsection{CT imaging geometry}

We consider a typical  
2D fan-beam geometry with $360^\circ$ circular source path and $\nv$ equi-angular projections, or views, onto a curved detector.
We consider a square domain of $\imside \times \imside$ pixels, and
due to rotational symmetry we restrict the region-of-interest to be within a disk-shaped mask
inside the square domain consisting of approximately
$\imsize = \lceil \pi/4 \cdot \imside^2 \rceil$ pixels.
The source-to-center distance is $2\imside$, and the fan-angle $28.07^\circ$ is set to precisely illuminate the disk-shaped mask.
The detector consists of $2\imside$ bins,
so the total number of measurements is $\sinosize = 2\imside\nv$.
The $\sinosize \times \imsize$ system matrix $A$ is
computed by means of the MATLAB package AIR Tools~\cite{Hansen2012}.

\subsection{Sparse image classes}

By an \emph{image class} we mean a set of test images described by a set of specifications, such that we can generate random realizations from the class.
We refer to such an image as an \emph{image instance} from the class, and multiple image instances from the same class form an \emph{image ensemble}.

\begin{figure}[htb]
\centering
\includegraphics[width=0.9\linewidth]{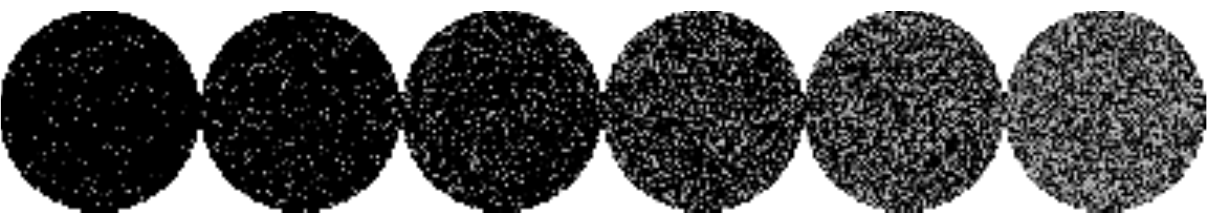}
 \includegraphics[width=0.9\linewidth]{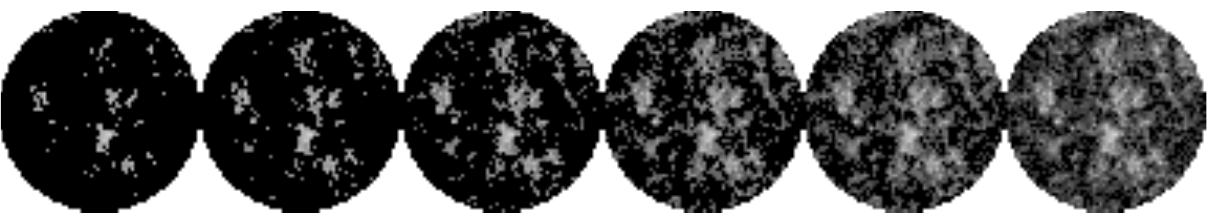}
 \includegraphics[width=0.9\linewidth]{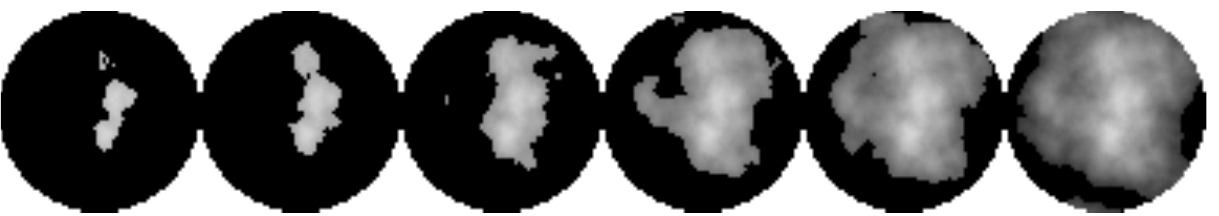}
\caption{Top, middle, bottom: Sparse image instances from \spikes{}, \onepower{}, \twopower{} classes. 
Left to right: $5\%$, $10\%$, $20\%$, $40\%$, $60\%$, $80\%$ non-zeros.
Gray-scale $[0,1]$.
\label{fig:spikes_instances}}
 \end{figure}

For the \spikes{} image class, given an image size $\imsize$ and a target sparsity $\sparsity$, we generate an image instance 
as follows:
starting from the zero image, 
randomly select $\sparsity$ pixel indices, and set each select pixel value to a random number from a uniform distribution over $[0,1]$.
Figure \ref{fig:spikes_instances} shows examples of \spikes{} image instances for
varying $\relsparsity$.
This class is deliberately designed to be as simple as possible and it does not model any particular application;
it is solely used to study the generic case of having a sparse image. The \spikes{} class, which is known as ``random-valued impulse noise'' in the signal processing literature, is often considered in CS studies, e.g., \cite{DonohoTanner:2009}, and therefore studying it in the setting of CT is natural.

The \ppower{} class is a more realistic image class that models background tissue in the female breast. 
The idea is to introduce structure to the pattern of non-zero pixels by creating
correlation between neighboring pixels.
Our procedure is based on \cite{Reiser2010} followed by thresholding to obtain many zero-valued pixels. 
The amount of structure is governed by a parameter~$p$\/:
\begin{enumerate}
 \item Create an $\imside\times\imside$ phase image $P$ with values drawn from a Gaussian
  distribution with zero mean and unit standard deviation.
 \item Create an $\imside\times\imside$ amplitude image $U$ with pixel values
\begin{align*}
 U(i,j) &= \left( \left(\frac{2i-1}{N}-1\right)^2 + \left(\frac{2j-1}{N}-1\right)^2 \right)^{-p/2}, \quad i,j = 1,\dots,\imside .
\end{align*}
 \item For all pixels $(i,j)$ compute $F(i,j) = U(i,j)e^{2\pi \hat{\imath} P(i,j)}$, with 
 $\hat{\imath} = \sqrt{-1}$.
 \item Compute square image as the magnitude of the DC-centered 2D inverse discrete Fourier transform
  of~$F$.
 \item Restrict the square image to the disk-shaped mask.
\item Keep the $\sparsity$ largest pixel values and set the rest to $0$.
\end{enumerate}
Figure \ref{fig:spikes_instances} shows examples of images from classes \onepower{} and \twopower{}.
Both have more structure than the \spikes{} images, and the structure
increases with $p$.

We do not claim that our image classes are fully realistic models. As our goal is to study how image sparsity affects the number of projections sufficient for recovery, we simply consider a selection of qualitatively different image classes. Development of more realistic image classes for specific applications is beyond the scope of the present work.

\subsection{Robust solution of optimization problem} \label{sec:reliable}

Given a numerically computed solution, the robustness of the decision regarding recovery depends on the accuracy of the solution.
False conclusions may result from incorrect or inaccurate solutions.
To robustly solve the optimization problem \ellone{} we must therefore
use a numerical method which gives a reliable indication of whether a correct solution,
within a given accuracy, has been computed.
Our choice is the package MOSEK \cite{MOSEK}, which uses a primal-dual interior-point method and
issues warnings in the rare case that an accurate solution can not be computed. For all problem instances in our studies, MOSEK returned a certified accurate solution.

To solve \ellone{} using MOSEK we recast it as the linear program 
\begin{equation}
 \argmin_{\im,w}\; \mathbf{1}^T w \quad \text{s.t.} \quad \sysmat \, \im = \sino \quad \text{and} \quad -w \leq \im \leq w,
\end{equation}
where $\mathbf{1} \in \mathbb{R}^\imsize$ and the inequality constraints imply non-negativity of $w$.

\subsection{Simulations}\label{sec:simulations}

Using the presented imaging model, method for generating sparse test images,
and robust optimization algorithm we carry out
simulation studies of recoverability within an image class.
A single recovery simulation consists of the following steps, assuming image sparsity $k$ and number of projections $\nv$:
\begin{enumerate}
 \item Generate $\sparsity$-sparse test image $\imorig$,
 \item generate system matrix $A$ using $\nv$ projections,
 \item compute perfect data $\sino = A\imorig$,
 \item solve \ellone{} numerically to obtain $\imopt$, and
 \item test for recovery numerically using 
\begin{align}
 \frac{\normtwo{\imopt - \imorig}}{\normtwo{\imorig}} < \numthr,
 \label{eq:teststrongrecov}
\end{align}
\end{enumerate}
where the threshold $\numthr$ is chosen based on the chosen accuracy of the optimization algorithm; empirically we found $\numthr = 10^{-4}$ to be well-suited in our set-up. 

We wish to study how the number of projections sufficient for \ellone{} recovery depends on the image sparsity. In order to make comparisons across image size we introduce the following normalized measures of sparsity and sampling. 
For a given problem instance, we define the \emph{relative sparsity} as
\begin{equation} \label{eq:relsparsity}
 \hbox{\bf relative sparsity:} \qquad \relsparsity = \sparsity/\imsize .
\end{equation}
In our studies we will generate test images of a desired relative sparsity $\relsparsity$ and set the sparsity as $k = \mathrm{round}(\relsparsity\imsize)$.
We let the \emph{sufficient projection number} $\nvsuf$ denote the smallest number
of projections that causes $\sysmat$ to have full column rank. At $\nv \geq \nvsuf$ the original image is the unique feasible point and hence minimizer no matter what objective is minimized and we therefore use $\nvsuf$ as a reference point of \emph{full sampling}. 
For a given problem instance, the \emph{number of projections for \ellone{} recovery} $\nrecov$ denotes the smallest number of projections
for which 
recovery is observed for all $\nv \geq \nrecov$.
We define the
\begin{eqnarray}
 \hbox{\bf relative sampling:} & \qquad & \relsampling\phantom{^{\ellone}} = \nv / \nvsuf , \\
 \hbox{\bf relative sampling for \ellone{} recovery:} & \qquad & \relsampling^{\ellone} = \nv^\ellone{} / \nvsuf.
\end{eqnarray}
For the test problems considered here, existence of an \ellone{} solution is guaranteed by the way we generate data. As mentioned in 
Section~\ref{sec:guarantees}, uniqueness is not guaranteed and for a given problem instance it cannot be known in advance whether the solution is unique.
The computed solution depends on the optimization algorithm, and therefore our conclusions of recoverability by \ellone{} are, in principle, subject to our use of MOSEK. We do not specifically check for uniqueness; however, in the event of infinitely many solutions, it is unlikely that any optimization algorithm will select precisely the original image, so we believe that our observations of recoverability based on solving \ellone{} correspond to existence of a unique solution.

\section{Simulation results}
\label{sec:numericalresults}

In this section we present our numerical results. We first establish that \ellone{} indeed can recover the original image from fewer than $\nvsuf$ projections. We then systematically study how $\nrecov$ depends on the image sparsity, image size, image class and finally the robustness to noise.

\subsection{Recovery from undersampled data}

\begin{figure}[htb]
\centering

\hspace*{0.5em}
\includegraphics[width=0.15\linewidth]{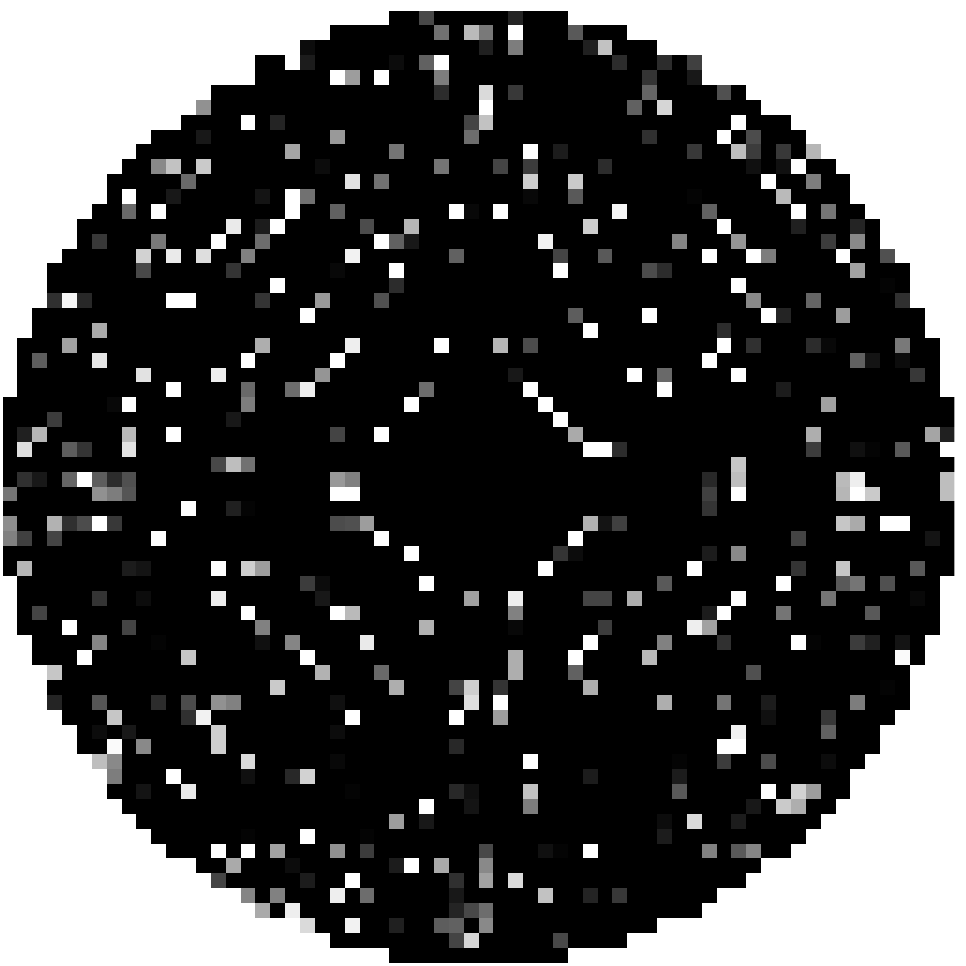}
\includegraphics[width=0.15\linewidth]{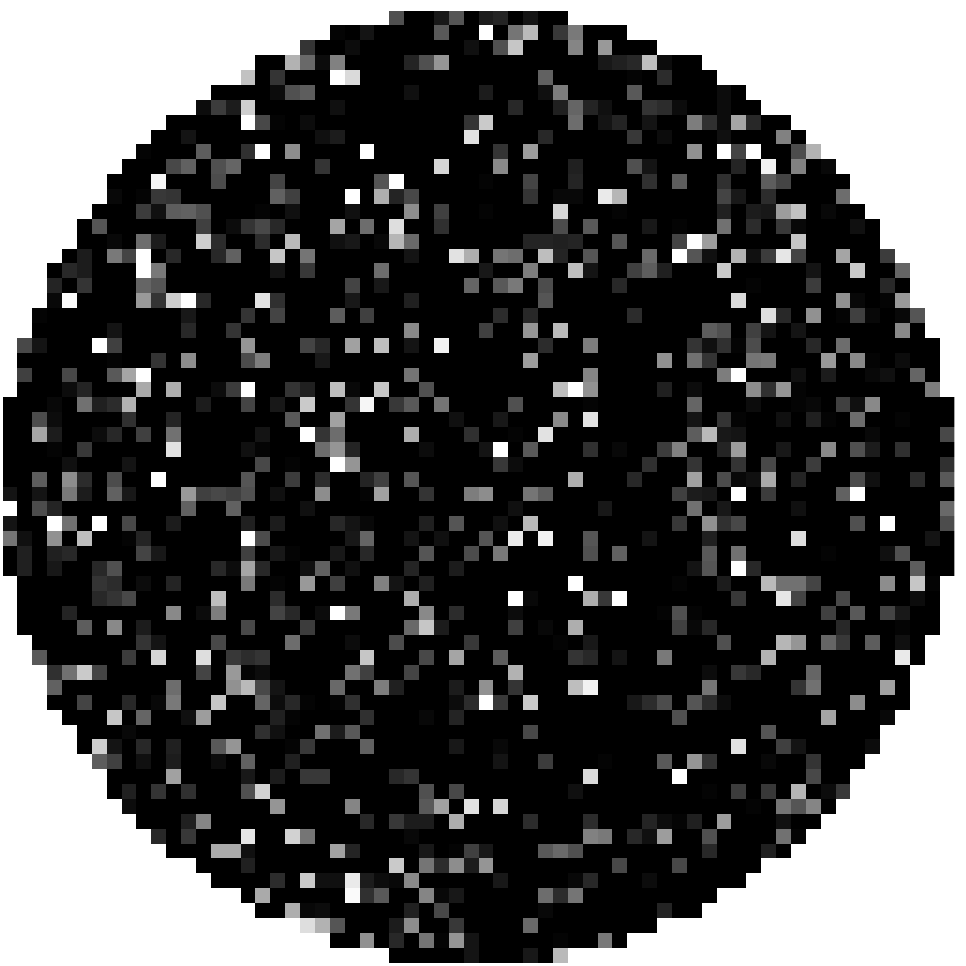}
\includegraphics[width=0.15\linewidth]{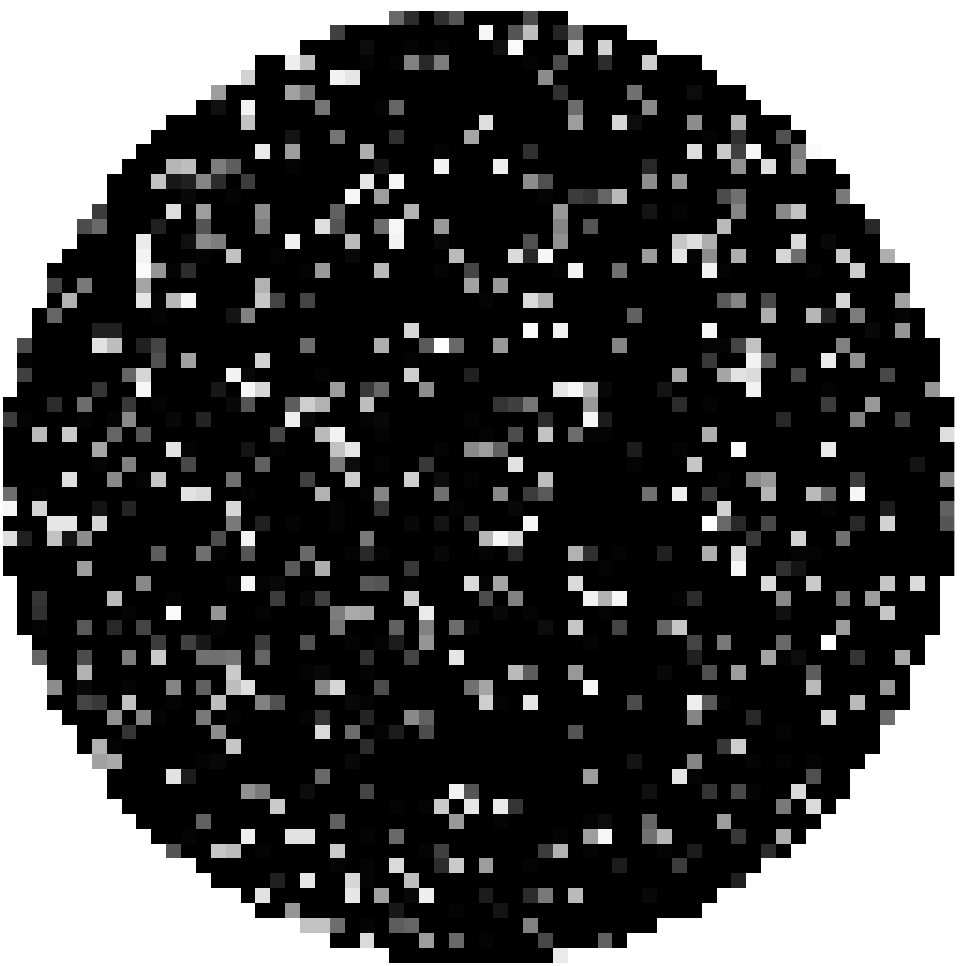}
\includegraphics[width=0.15\linewidth]{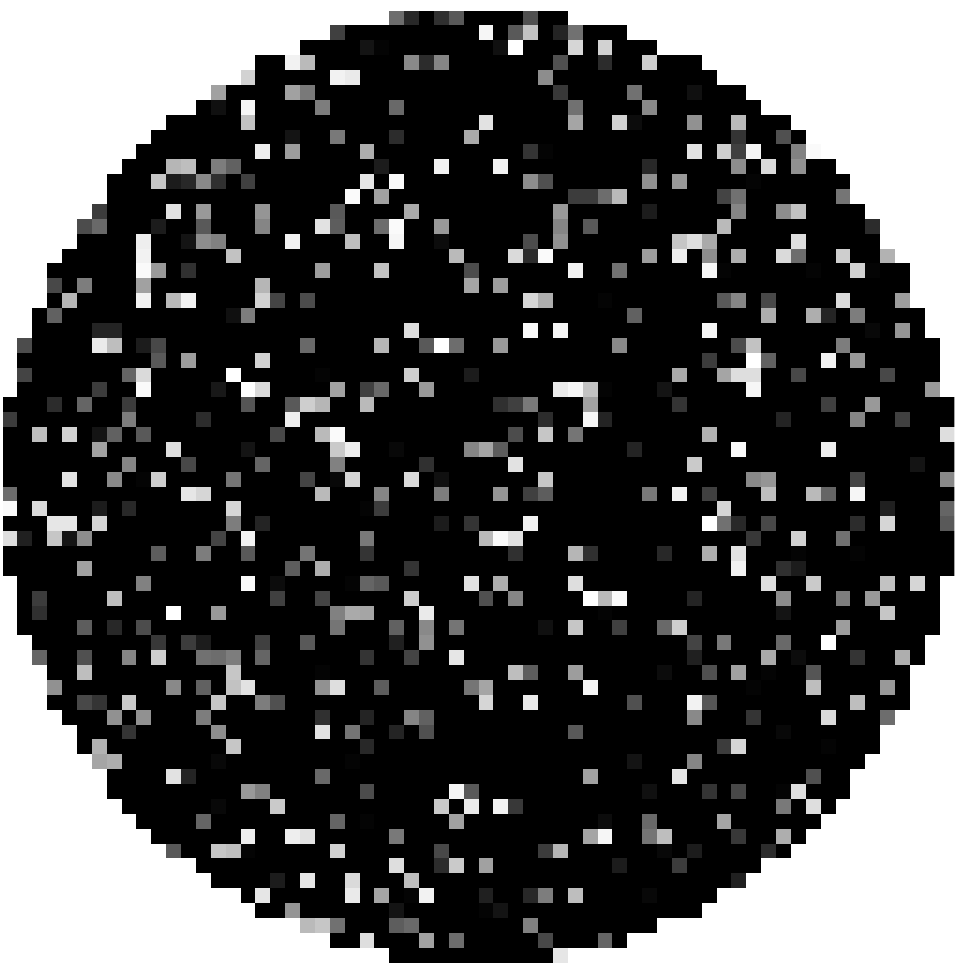}
\includegraphics[width=0.15\linewidth]{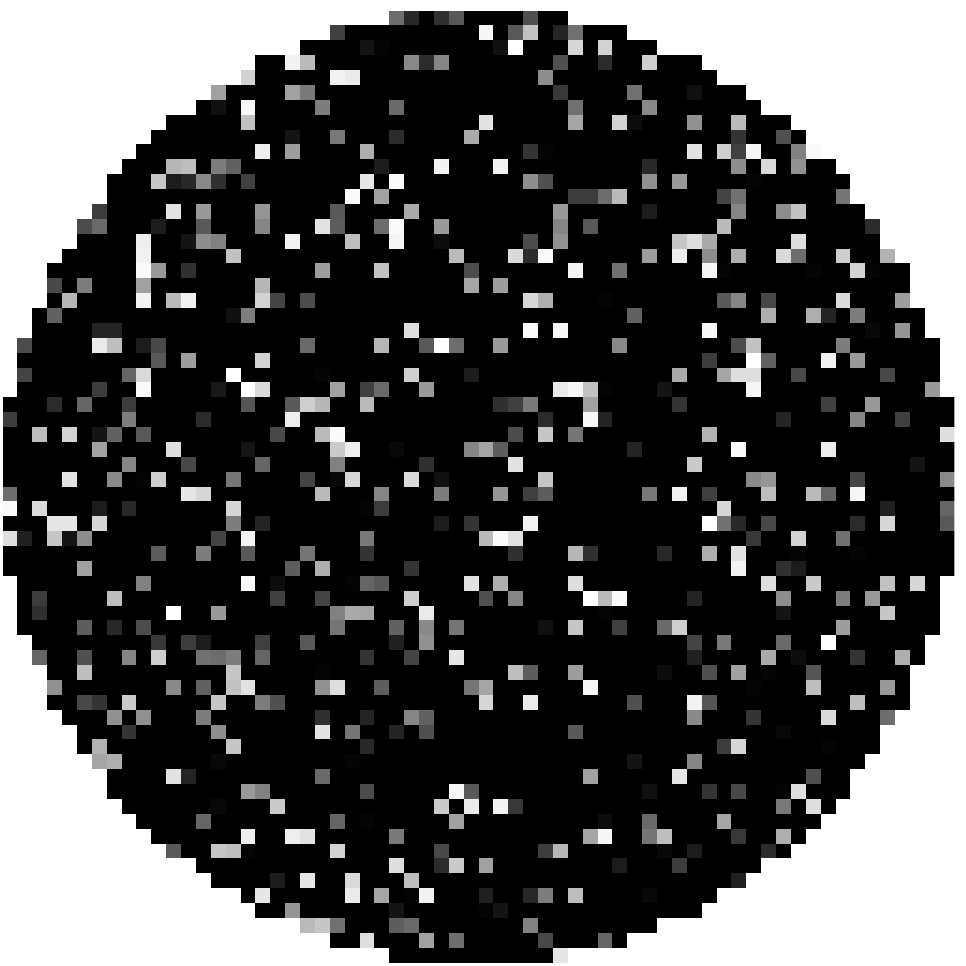}
\includegraphics[width=0.15\linewidth]{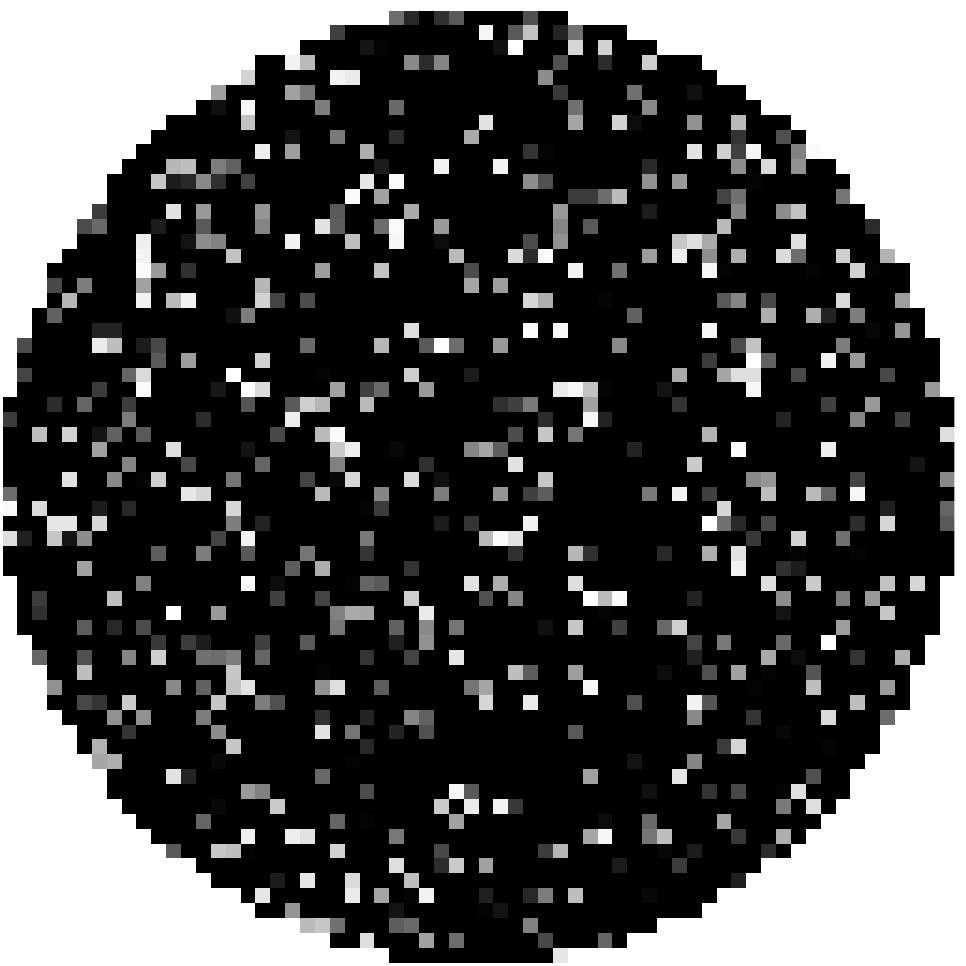}\\
\hspace*{0.5em}
\includegraphics[width=0.15\linewidth]{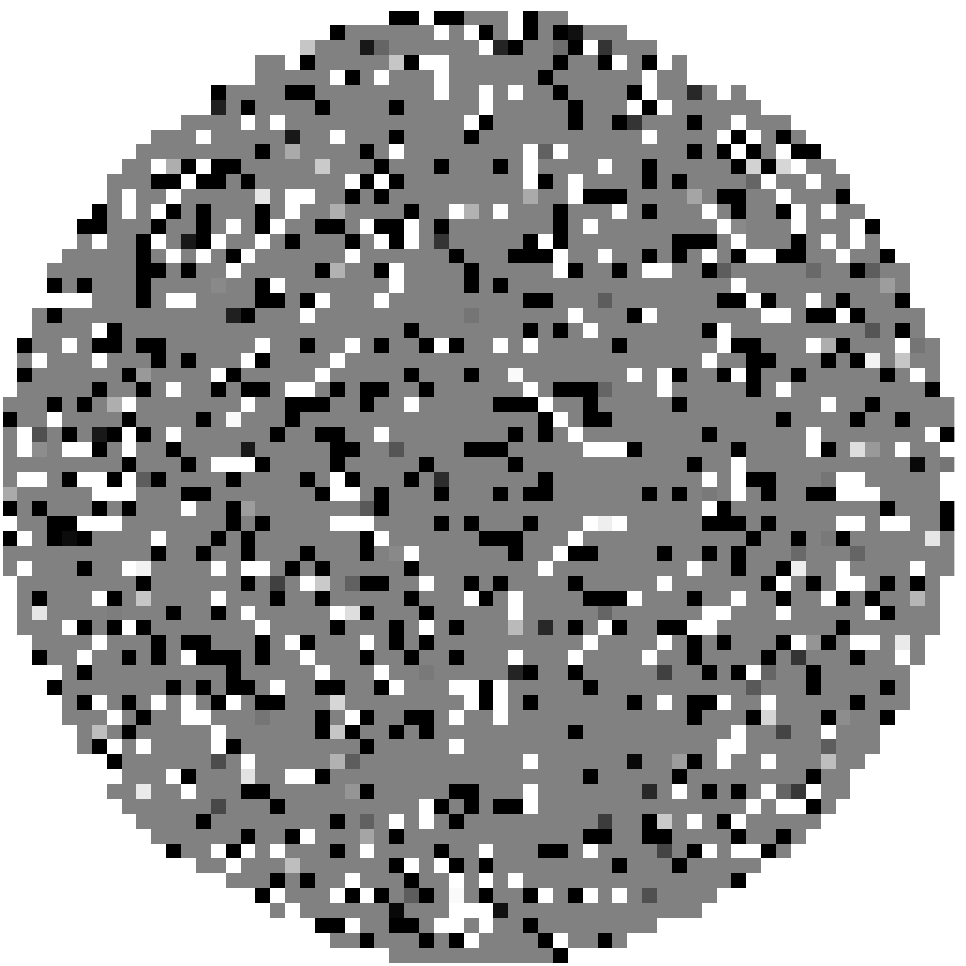}
\includegraphics[width=0.15\linewidth]{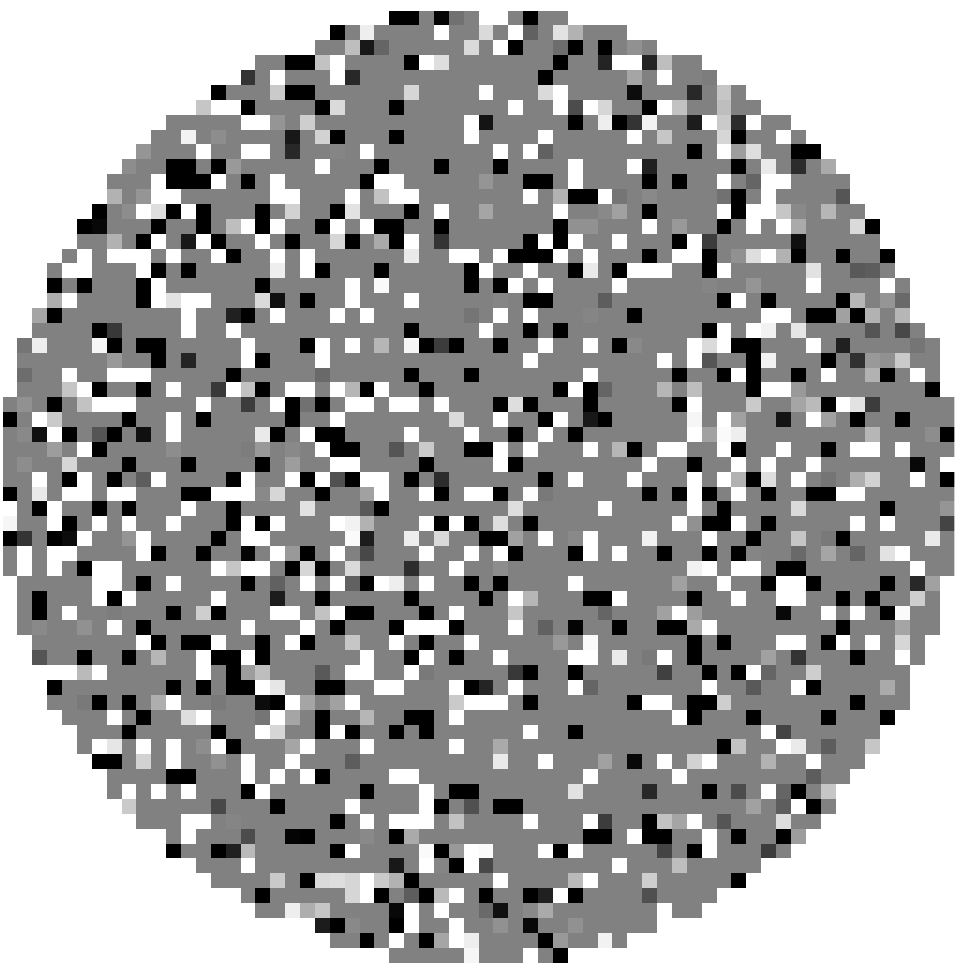}
\includegraphics[width=0.15\linewidth]{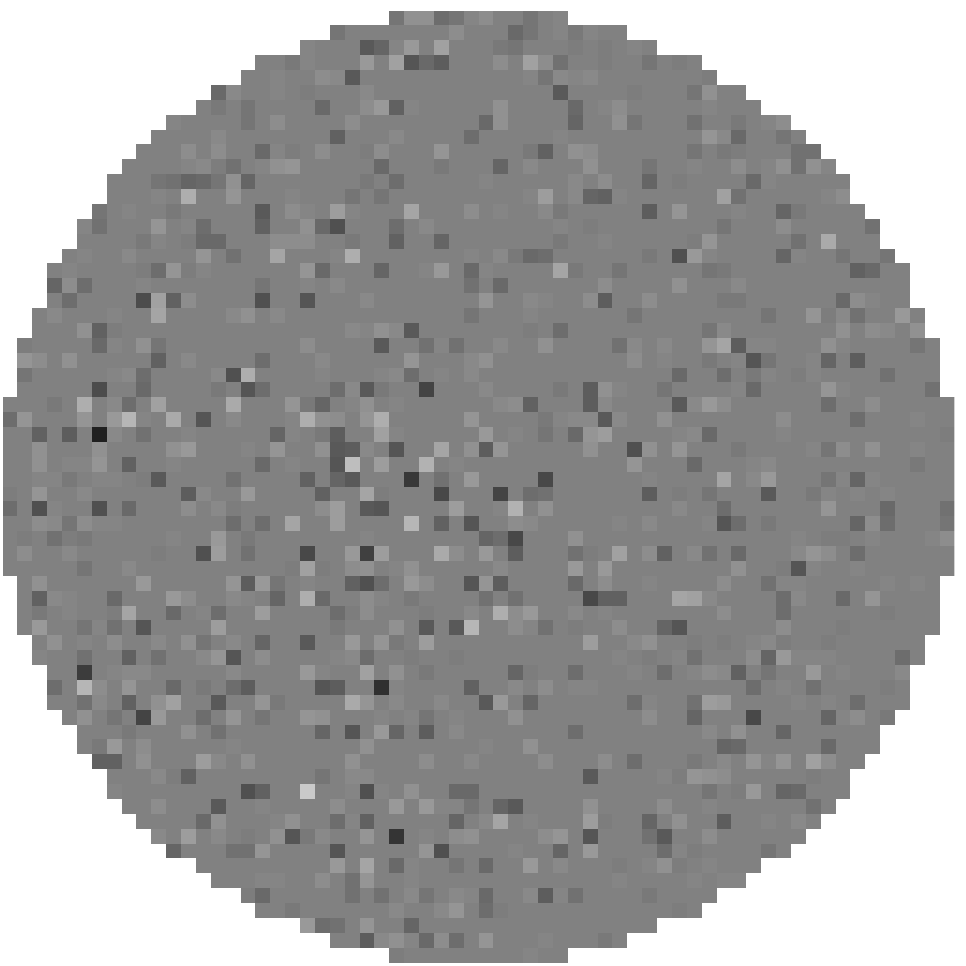}
\includegraphics[width=0.15\linewidth]{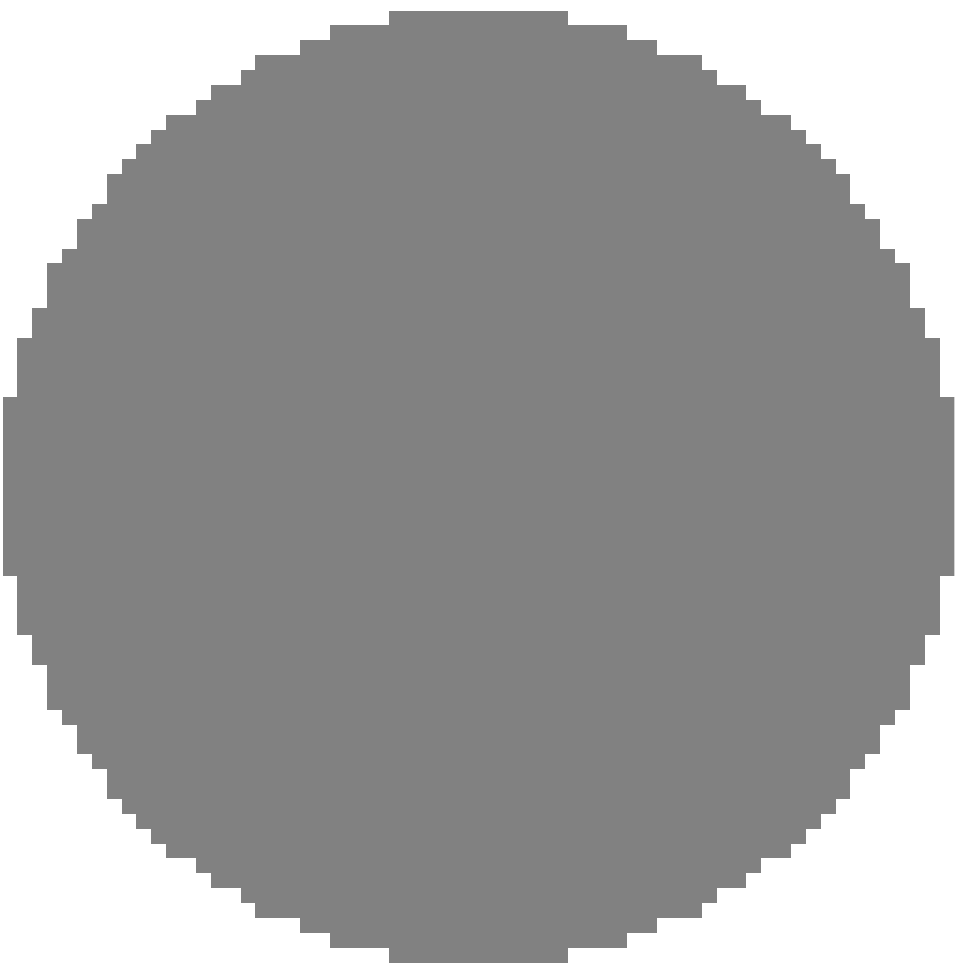}
\includegraphics[width=0.15\linewidth]{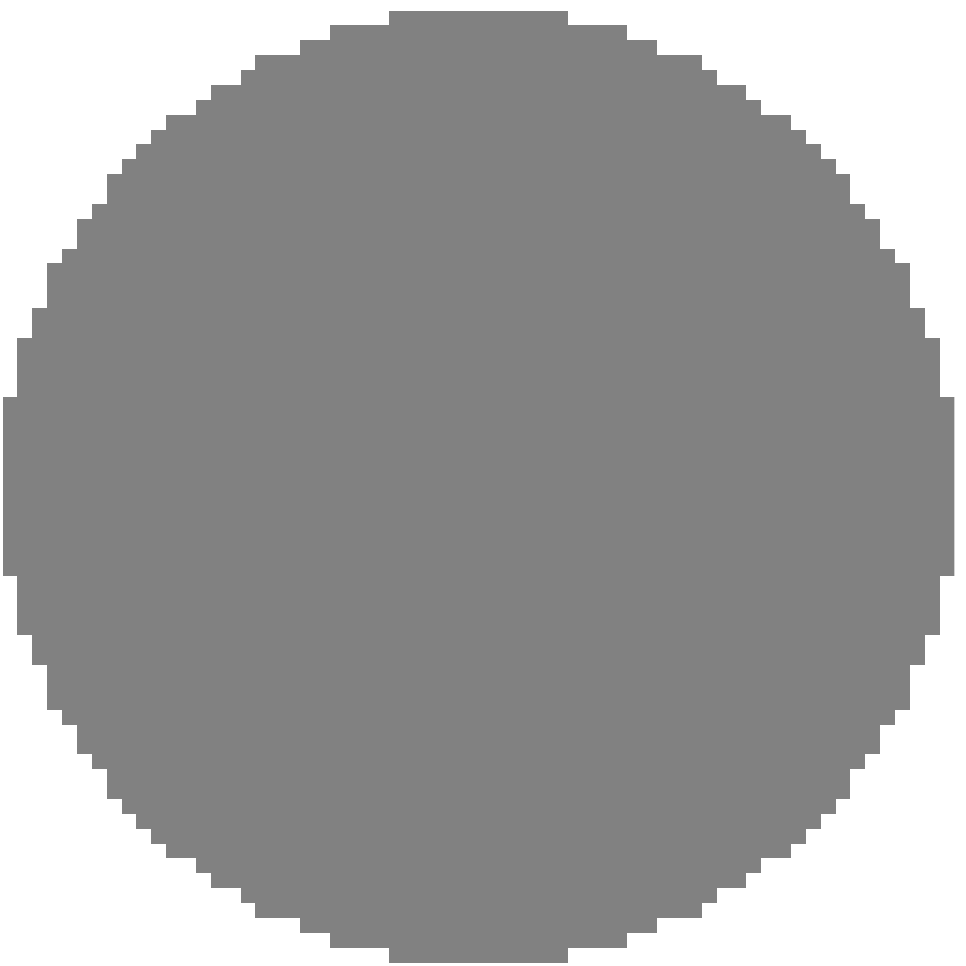}
\includegraphics[width=0.15\linewidth]{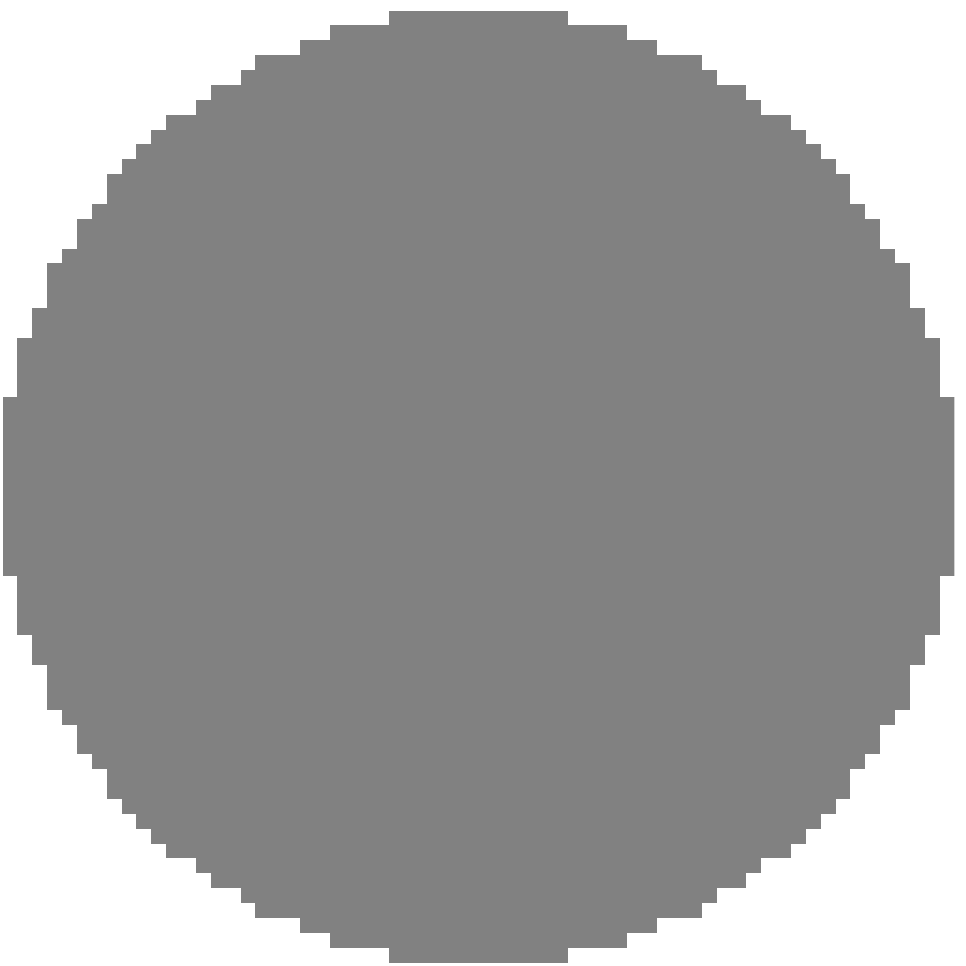}\\
\caption{\ellone{} reconstructions of \spikes{} image with $\imside=64$, relative sparsity $\relsparsity = 0.20$. 
Top row: \ellone{} reconstructions, gray-scale $[0,1]$. Bottom row: \ellone{} minus original image, gray-scale $[-0.1,0.1]$. Columns: $4$, $8$, $10$, $12$ ($=\nrecov$), $24$ and $26$ ($=\nvsuf$) projections.\label{fig:somereconstructions}}
 \end{figure}

To verify that \ellone{} is capable of recovering an image from undersampled CT measurements in our set-up, 
we use a \spikes{} image $\imorig$ with $\imside = 64$, leading to $N = 3228$
pixels in the disk-shaped mask.
The relative sparsity is set to $\relsparsity = 0.20$, which yields $646$ non-zeros. We consider reconstruction from data corresponding to $2,4,6,\ldots,32$ projections; the smallest and largest system matrices are of sizes $256\times{}3228$ and $4096\times{}3228$, respectively. 
At $\nv = 24$, the matrix is $3072\times{}3228$ and $\mathrm{rank}(A)=3052$, at $\nv = 25$ it is $3200\times{}3228$ and has rank $3185$, while at 
$\nv = 26$, the matrix is $3328\times{}3228$ and full-rank;
hence $\nvsuf = 26$.
Selected \ellone{} reconstructed images $\imopt$are shown in Figure~\ref{fig:somereconstructions} along with the error images $\imopt - \imorig$ to better visualize the abrupt drop in error when the image is recovered.
\ellone{} recovery occurs already at $\nv = \nrecov = 12$, where $\sysmat$ has size $1536\times{}3228$ and rank $1524$, i.e., a substantial undersampling relative to the full-sampling reference point of $\nvsuf = 26$.

 \begin{figure}[htb]
\centering
 \includegraphics[width=0.9\linewidth]{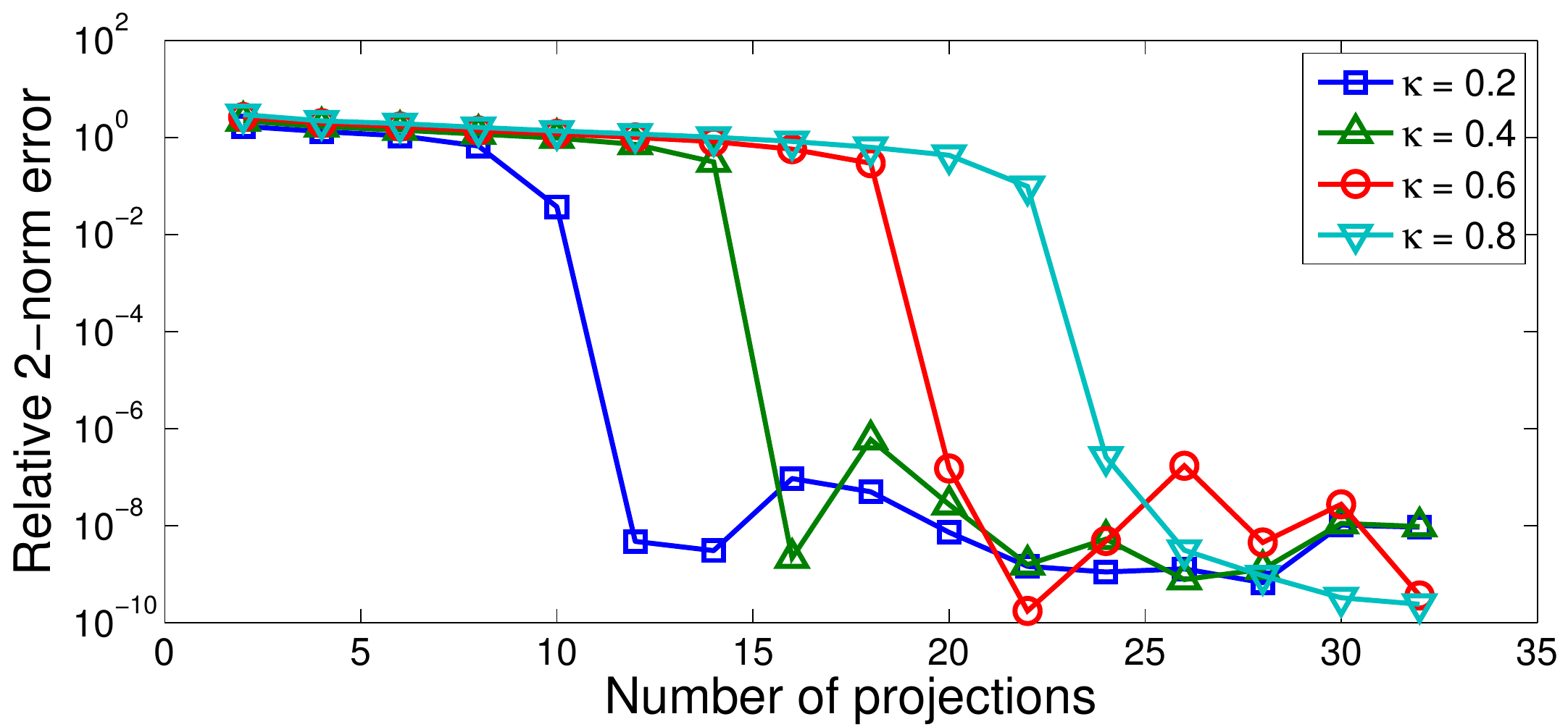}
\caption{Relative 2-norm error from \eqref{eq:teststrongrecov} of \ellone{} reconstructions vs.\ number of projections for \spikes{} images with relative sparsity values $\relsparsity = 0.2$, $0.4$, $0.6$, $0.8$. The $\kappa = 0.2$ image is recovered at $\nv = 12$ as also seen in Fig. \ref{fig:somereconstructions}. The relative errors of recovered images are non-zero due to the numerical accuracy chosen for MOSEK.
\label{fig:reconstruction_errors}}
 \end{figure}

To see how relative sparsity affects recovery,
we repeat the experiment for $\relsparsity = 0.4$, $0.6$ and $0.8$.
Figure~\ref{fig:reconstruction_errors} shows the relative 2-norm error from \eqref{eq:teststrongrecov} for the \ellone{} reconstructions 
as a function of $\nv$. In all cases an abrupt drop in error to a numerically accurate reconstruction is observed and 
the relative sampling for \ellone{} recovery $\nrecov$ changes from $\nrecov = 12$ at $\relsparsity = 0.2$ to $16$, $20$ and $24$.
This indicates a very simple relation between sparsity and number of projections for \ellone{} recovery.

\subsection{Recovery phase diagram}
\label{sec:phasediagram}

\begin{figure}[htb]
\centering
 \includegraphics[width=0.415\linewidth,clip,trim=0cm 0cm 0cm 1.7cm]{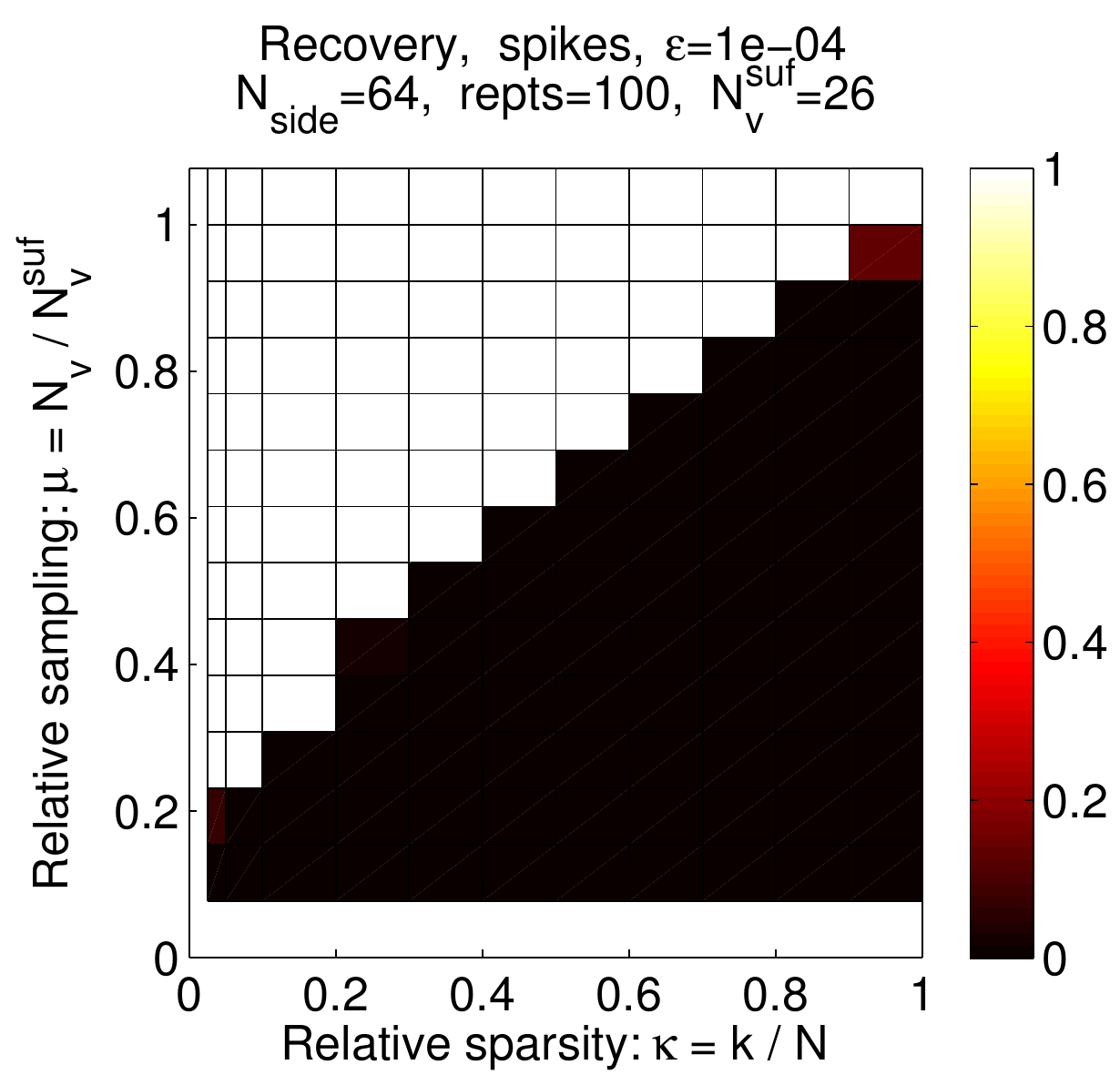}
\includegraphics[width=\wtwo\linewidth,clip,trim=0cm 0cm 0cm 2cm]{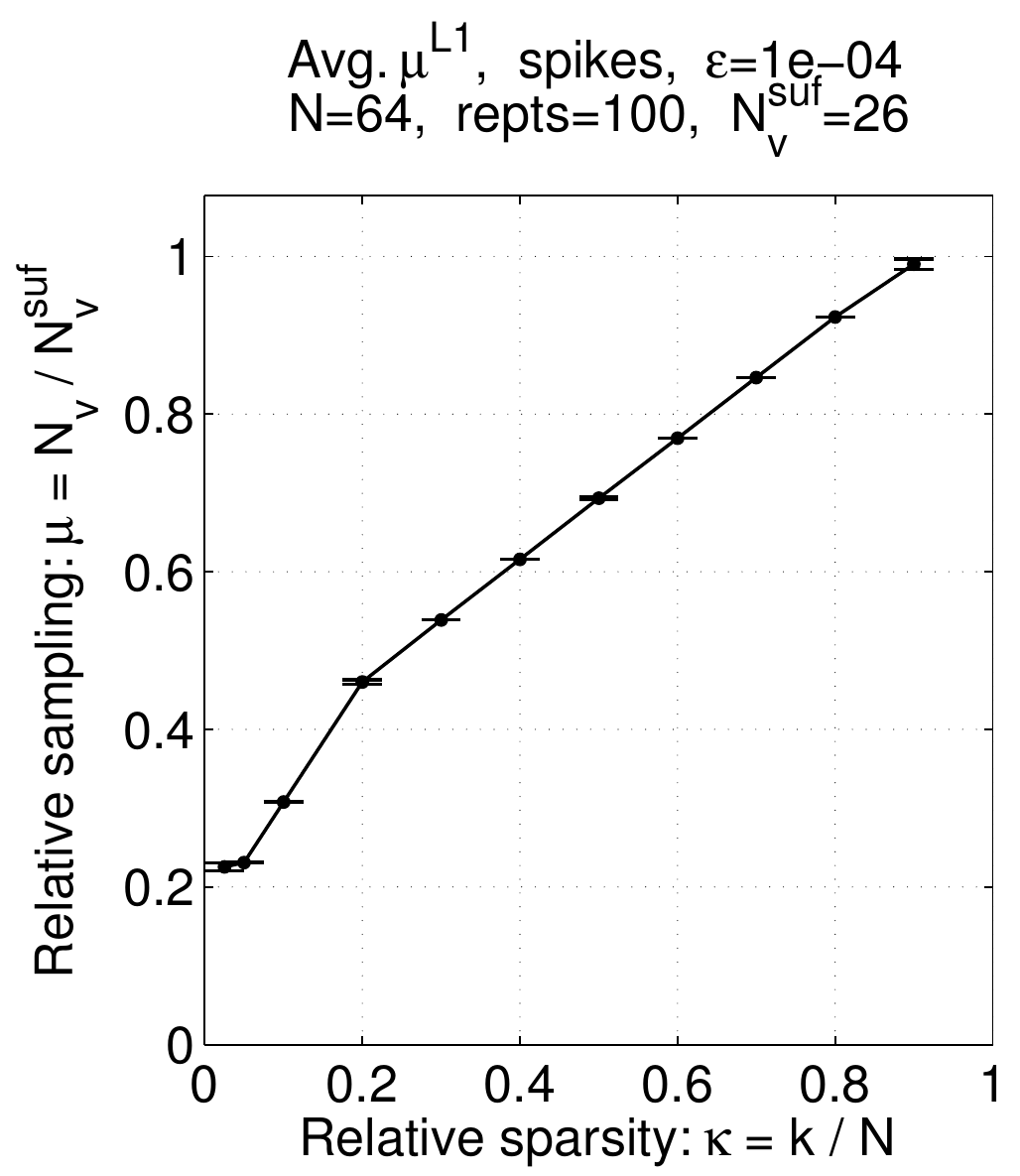}
\caption{Left:\ Phase diagram for \spikes{} images ($\imside=64$) showing percentage of image instances recovered by \ellone{} as function of relative sparsity and relative sampling. 
Right:\ Average relative sampling for \ellone{} recovery and its $99\%$ confidence interval. \label{fig:averagecaserecovery}}
 \end{figure}

In general, we can not expect
all image instances of the same relative sparsity to have the same $\nrecov$.
We can study the variation within an image class by determining the number of projections for \ellone{} recovery of an ensemble of images of different sparsity.
At each of the relative sparsity values $\relsparsity=0.025$, $0.05$, $0.1,0.2,\ldots,0.9$ we create $100$ image instances of the \spikes{} class. For each instance we compute the \ellone{} reconstruction from $\nv = 2, 4, \ldots, 32$ projections. Based on the relative $2$-norm error in \eqref{eq:teststrongrecov} we assess whether the original is recovered. Figure~\ref{fig:averagecaserecovery} (left) shows the percentage of recovered instances as function of relative sparsity $\relsparsity$ and relative sampling $\relsampling$. Each square corresponds to the $\relsparsity$ value at the left edge of the square and the $\relsampling$ value at the bottom edge. We refer to this plot as a \emph{phase diagram}.
The phase diagram shows two distinct regions: the lower right black one, in which no image instances were recovered, and the upper left white one, in which all images were recovered.

An important observation is the sharp phase transition from non-recovery to recovery, meaning that the number of projections for \ellone{} recovery is essentially constant for same-sparsity images of the \spikes{} class. A sharp phase transition is often seen in CS \cite{DonohoTanner:2009}, but to the best of our 
knowledge has not been reported for CT-matrices before, and we therefore believe the sharp transition to be a novel observation. The sharpness of the transition is perhaps better appreciated in Figure~\ref{fig:averagecaserecovery} (right), 
which shows the \emph{average} $\relsampling^\ellone{}$ over all instances at each $\relsparsity$ and its $99\%$ confidence interval estimated using the empirical standard deviation, illustrated by errorbars. The confidence intervals are very narrow, in fact, in several cases of width zero, due to zero variation of $\nrecov$, which agrees with the visual observation of a sharp transition. 

The relative sampling for recovery $\relsampling^\ellone{}$ increases monotonically with the relative sparsity $\relsparsity$. %
As $\relsparsity \rightarrow 0$, also $\relsampling^\ellone{} \rightarrow 0$, showing that extremely sparse images can be recovered from very few projections, i.e., highly undersampled data.
As $\relsparsity \rightarrow 1$, $\relsampling^\ellone{} \rightarrow 1$, confirming that \ellone{} does not admit undersampling for non-sparse signals. 
Furthermore, the phase diagram makes these observations \emph{quantitative}. 
Assume, for example, that we are given an image of relative sparsity $\relsparsity = 0.1$,
how many projections would suffice for recovery?
The phase diagram 
shows that at $\relsparsity = 0.1$ the average
$\relsampling^\ellone{} = 0.31$, corresponding to $\nrecov =8$ projections. 
Or conversely, the maximal relative sparsity that, on average, allows recovery from $8$ projections is $\kappa = 0.1$.

\begin{figure}[htb]
\centering
 \includegraphics[width=0.45\linewidth,clip,trim=0cm 0cm 0cm 2cm]{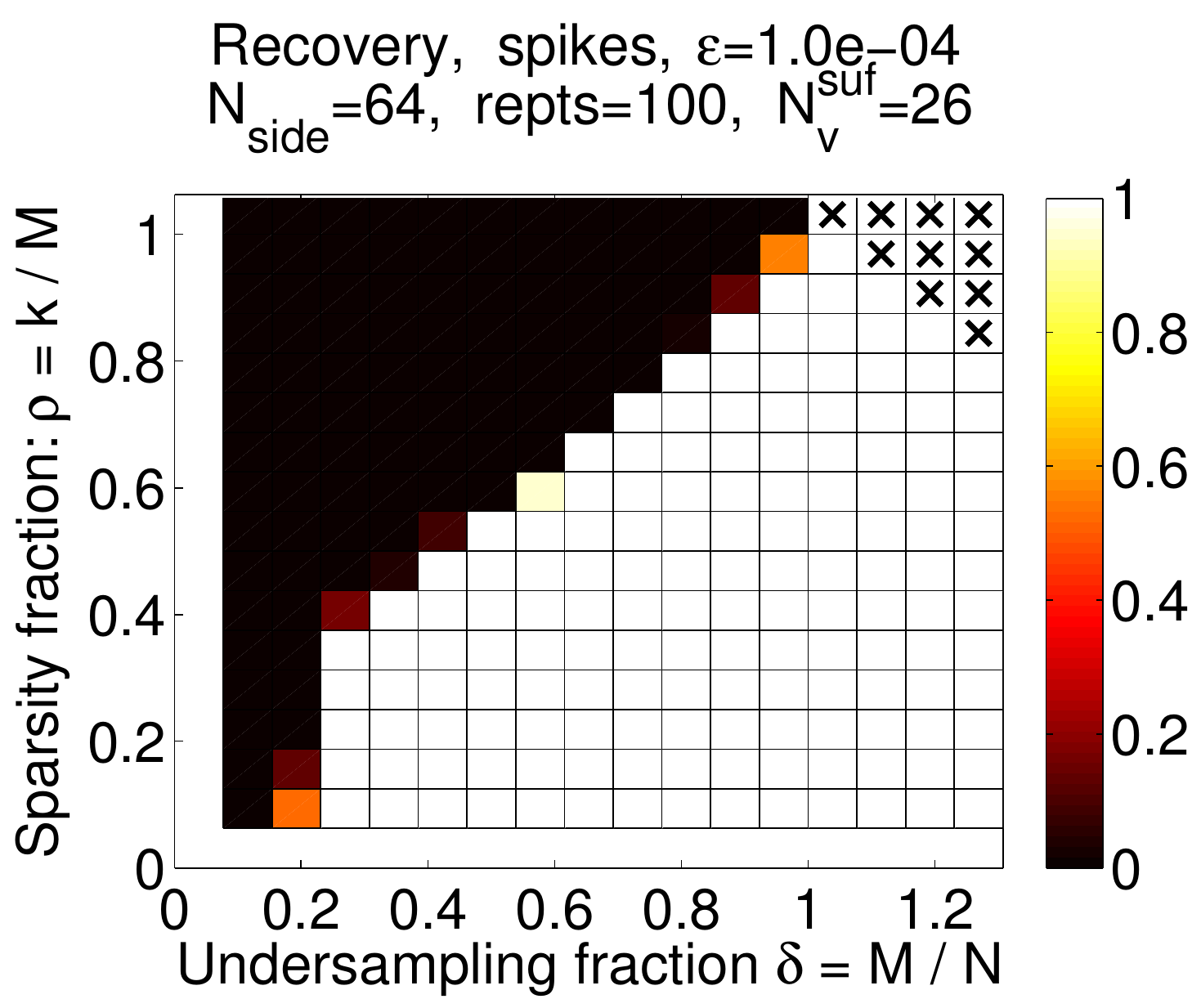}
\caption{Donoho-Tanner phase diagram
for \spikes{} images.
\label{fig:dt_diagram}}
 \end{figure}

We note that the phase diagram introduced by Donoho and Tanner \cite{DonohoTanner:2009} is slightly different from the one presented here. 
The DT diagram is parametrized by the sparsity fraction $\rho = \sparsity / M$, i.e. normalized by the number of measurements, not pixels, and undersampling fraction $\delta = M/\imsize$. 
We create a DT phase diagram for our set-up by assessing recovery for $100$ instances at $\nv = 2$, $4$, $\ldots$, $32$ and $\rho = 1/16, 2/16, \ldots, 16/16$, see Figure~\ref{fig:dt_diagram}. The DT phase diagram confirms the observations from Figure~\ref{fig:averagecaserecovery}, in particular the sharp phase transition.

We find that phase diagrams such as the one in Figure~\ref{fig:averagecaserecovery} are more intuitive to interpret because the quantities of interest, namely sparsity and sampling, can be read more directly of the axes. Furthermore, normalizing the sparsity by the number of samples as in the Donoho-Tanner phase diagrams leads to two minor issues: First, where in Figure~\ref{fig:averagecaserecovery} each column is based on $100$ specific instances of a certain sparsity, each square in the Donoho-Tanner diagram contains $100$ new instances, and hence sufficient sampling for specific instances is not addressed. Second, having $M\geq \imsize$ is common in CT, for example if using a large number of projections compared to the number of pixels, but as $\sparsity$ cannot be larger than $\imsize$, the upper right square cannot be realized, so they are marked by crosses. In summary, we consider phase diagrams of the former type more convenient in the setting of CT, and for the remainder of the paper we will only 
show this 
type of phase diagram.

\subsection{Dependence on image size}\label{sec:imagesize}

\begin{figure}[htb]
\centering
 \includegraphics[width=\wone\linewidth]{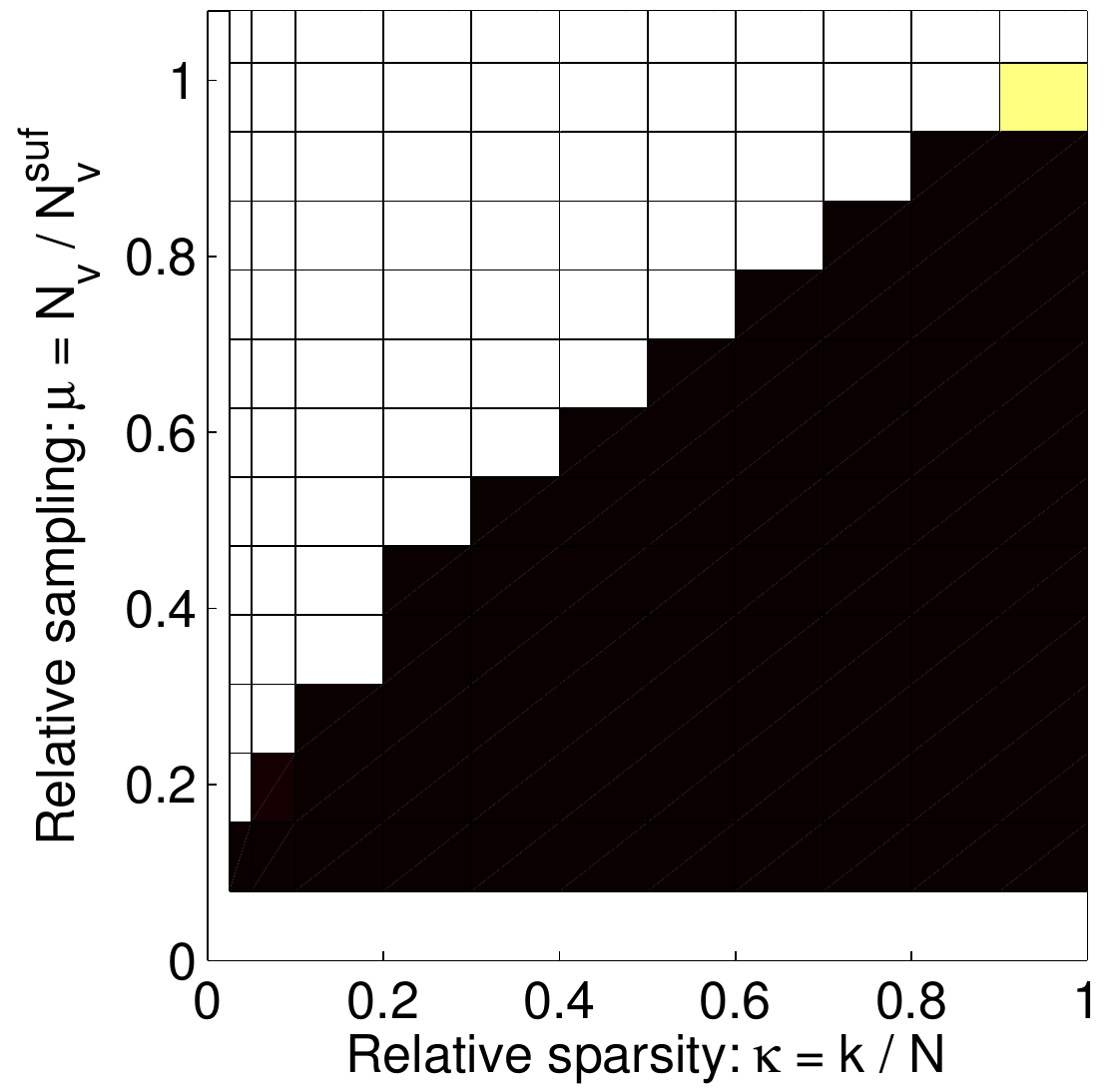}
.\includegraphics[width=\wone\linewidth]{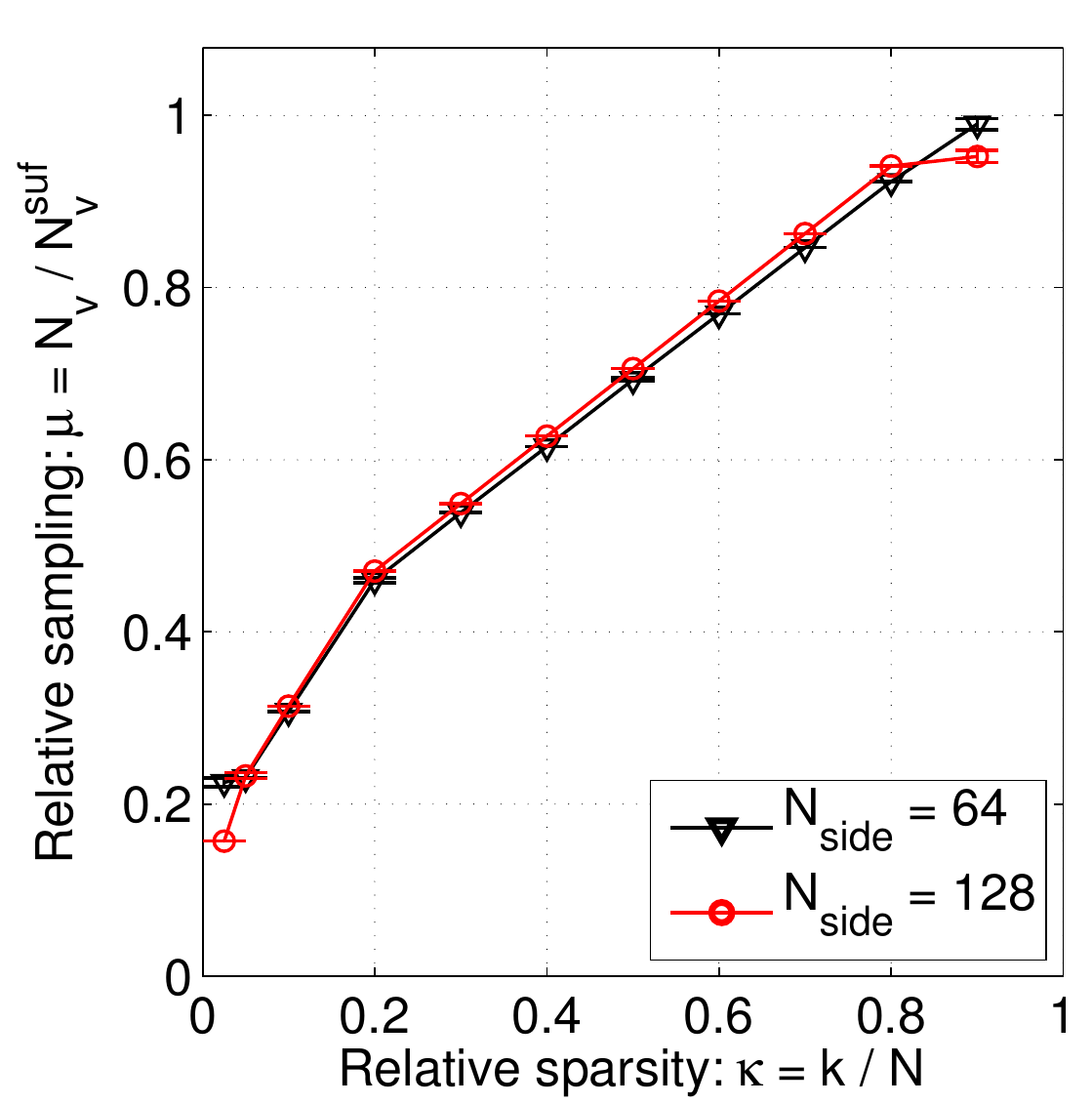} 
\caption{Phase diagram dependence on image size. Left:\ $\imside=128$. 
Right:\ Average relative sampling for \ellone{} recovery and $99 \%$ confidence intervals.\label{fig:RSS_imagesize}}
\end{figure}

To study how recoverability depends on image size,
we construct the phase diagram for $\imside = 128$, see Figure~\ref{fig:RSS_imagesize}.
For $\imside = 128$ we have $\nvsuf = 51$, so by taking $\nv = 4$,
$8$, \dots, $64$ we obtain approximately the same relative sampling values as for $\imside = 64$.

Overall, for $\imside = 128$ we see the same monotone increase in $\relsampling^\ellone{}$
with increasing~$\relsparsity$ as we did for $\imside = 64$.
The only differences are a generally sharper transition, i.e., narrower confidence intervals, as well as slightly better recovery at the extreme $\relsparsity$-values.
We conclude that, with appropriate normalization, the observed relation between the average number of projections for \ellone{} recovery and the image sparsity does not depends on the image size. Moreover, $\imside = 64$ is sufficiently large
to give
representative results that can be extrapolated to predict the sparsity-sampling relation for larger images.

\subsection{Dependence on image class}

\begin{figure}[htb]
\centering
 \includegraphics[width=\wone\linewidth]{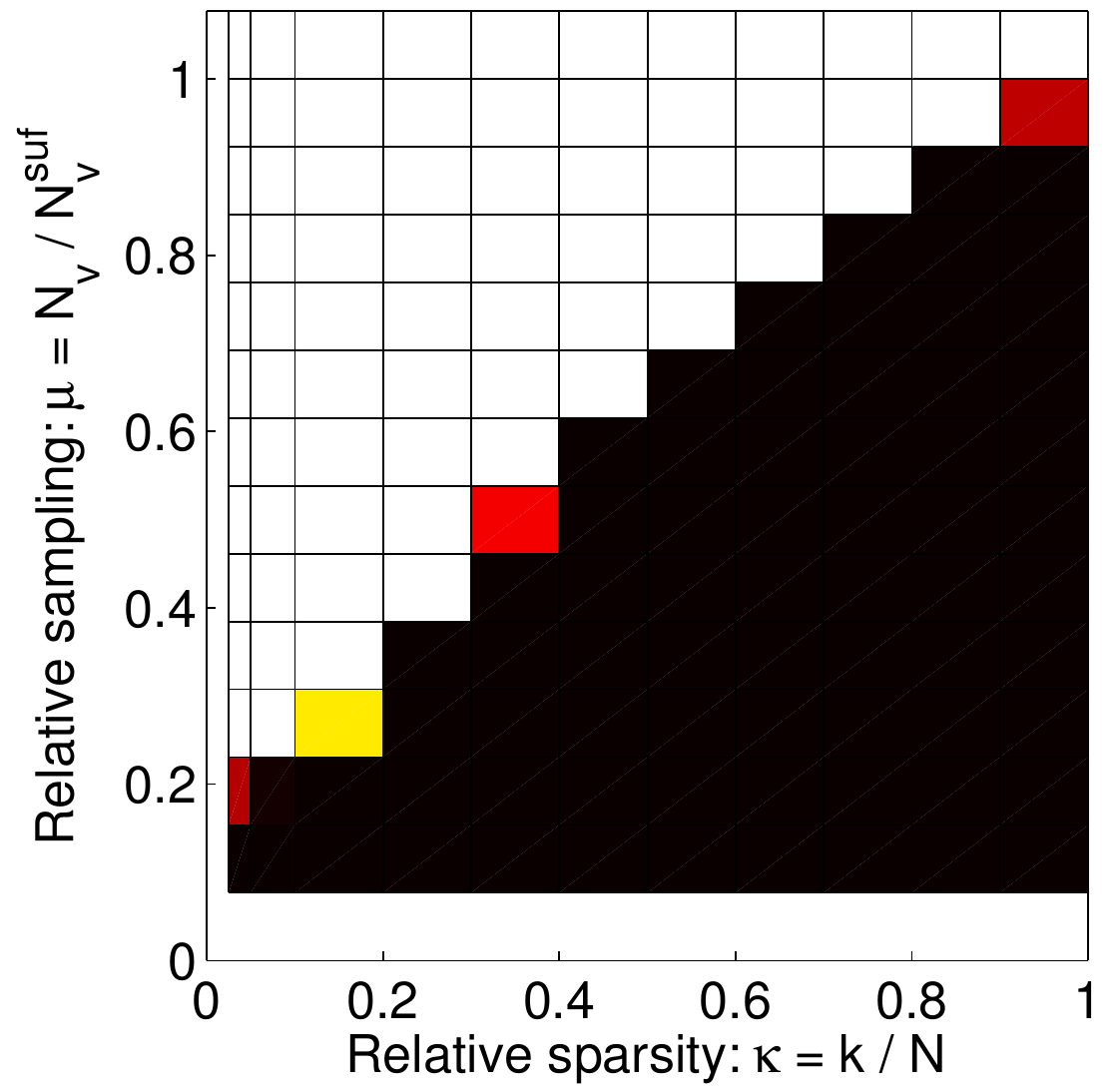}
 \includegraphics[width=\wone\linewidth]{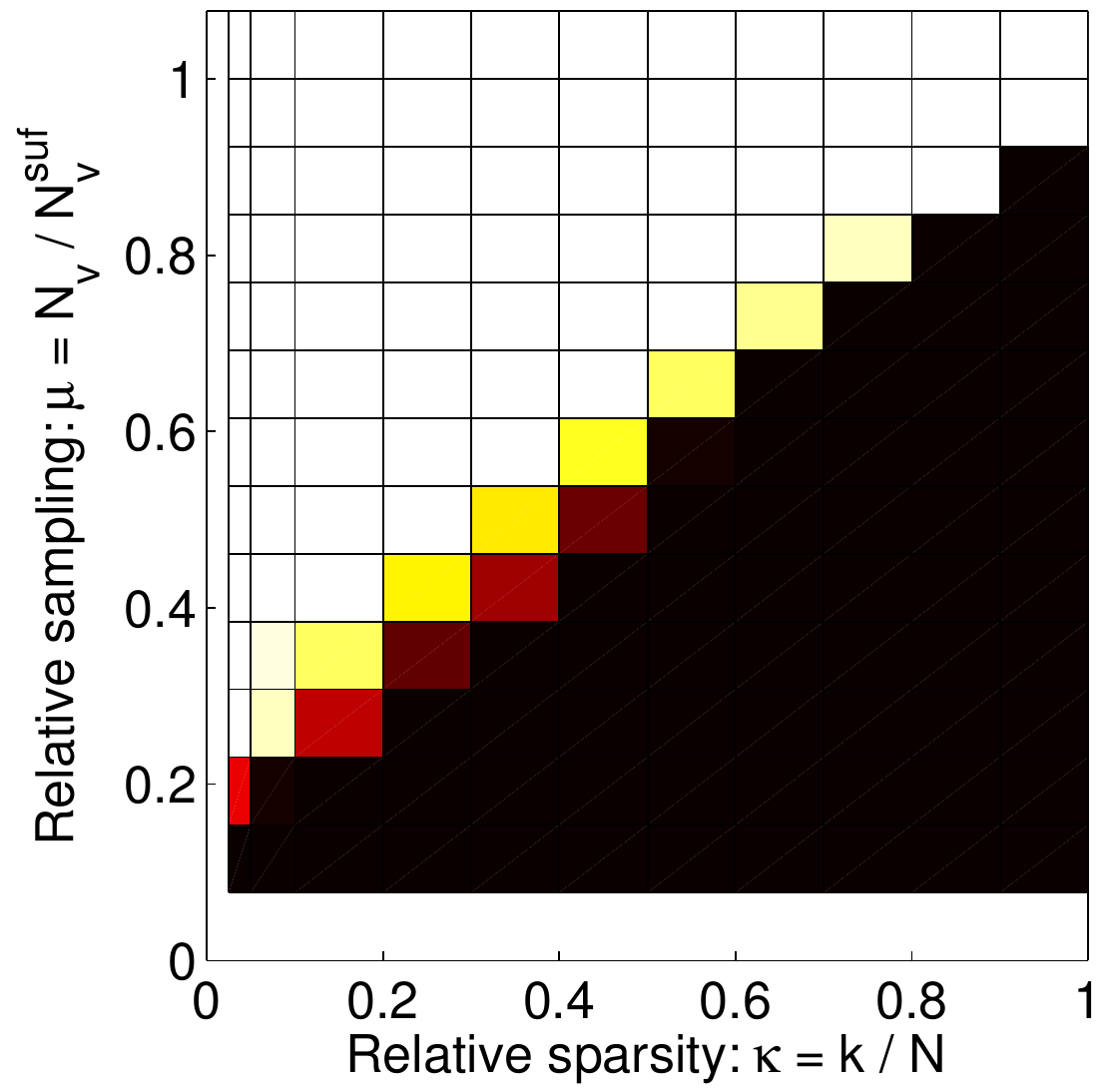}
\includegraphics[width=\wone\linewidth]{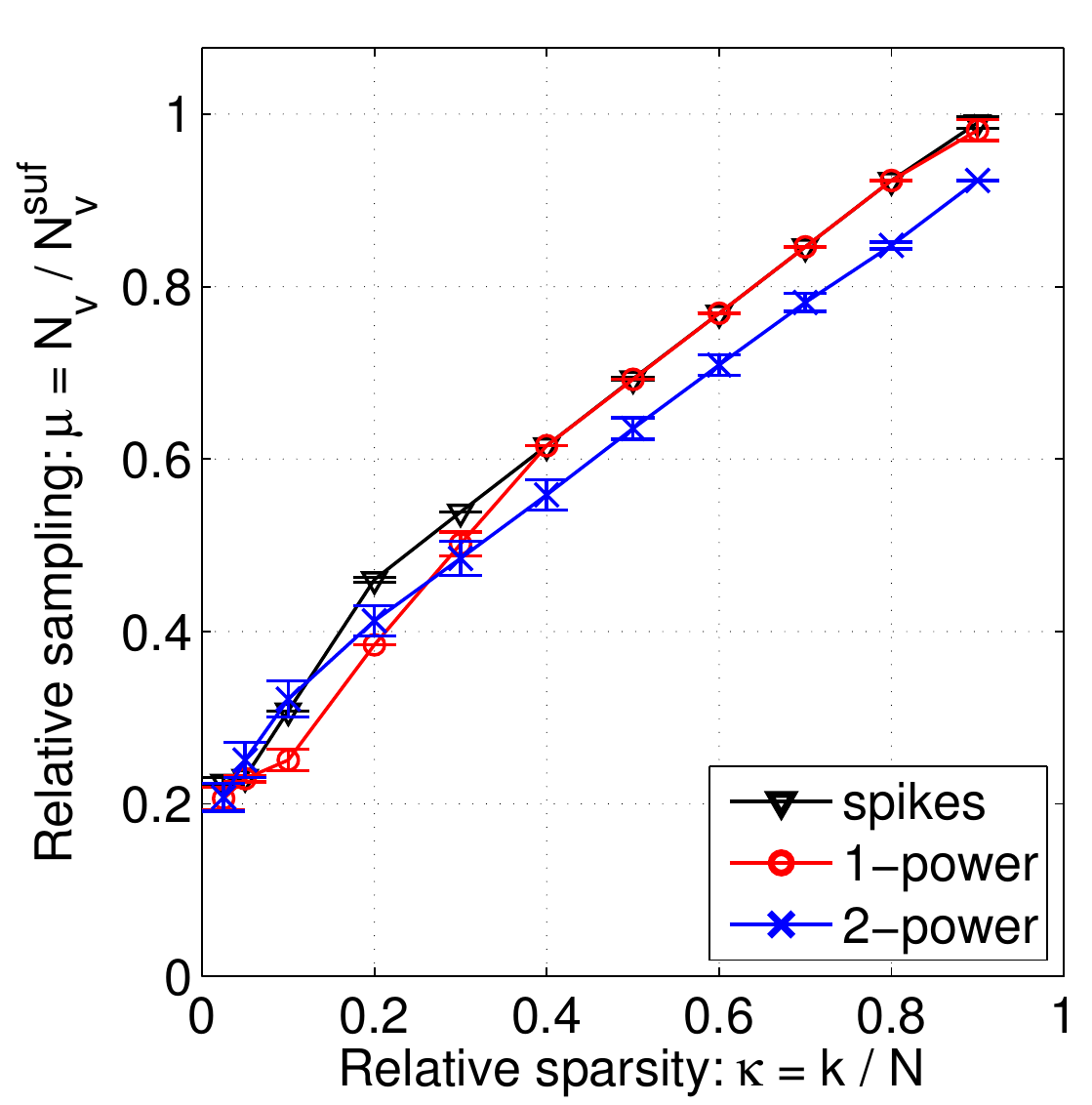}
\caption{Phase diagram dependence on image class. Left: \onepower{}. Center: \twopower{}. Right: Average sampling for \ellone{} recovery and $99\%$ confidence intervals. \label{fig:rss_ppower}}
 \end{figure}

Next, we study how the image class affects recoverability.
Figure \ref{fig:rss_ppower} shows phase diagrams for the \onepower{} and \twopower{}
classes for $\imside = 64$.
Comparing with \spikes{} in Figure \ref{fig:averagecaserecovery}
we observe similar overall trends but also some differences, which are more clearly seen in the plot of the average $\relsampling^\ellone{}$ for all classes in Figure~\ref{fig:rss_ppower} (right).
For \onepower{} the transition from non-recovery to recovery occurs at almost exactly the
same $(\relsparsity,\relsampling)$-values, except around $\kappa = 0.2$ where $\relsampling^\ellone{}$ is slightly lower, and is almost as sharp as for the \spikes{} class.
For \twopower{} the transition is more gradual, and occurs at lower
$\relsampling$-values for the mid and upper range of $\relsparsity$.

Based on these results we conclude that image classes with increasing structure on average admit recovery from a smaller number of projections but the
in-class recovery variability is larger. 
Thus, while recoverability is clearly tied to sparsity, the spatial correlation of the non-zero pixel locations also plays a role.

\section{Robustness to noise} \label{sec:noise}
In the previous section we have empirically established a relation between the number of projections for \ellone{} recovery and image sparsity in the noise-free setting.
A natural question
is whether and how the results generalize in the case of noisy data. 
Noise and inconsistencies in CT data are complex subjects arising from many different sources including scatter and preprocessing steps applied to the raw data before the reconstruction step. A comprehensive CT noise model is very application-specific and beyond our scope; rather we wish to investigate how recovery changes when subject to simple Gaussian white noise.

\begin{figure}[htb]
\centering
 \includegraphics[width=0.345\linewidth]{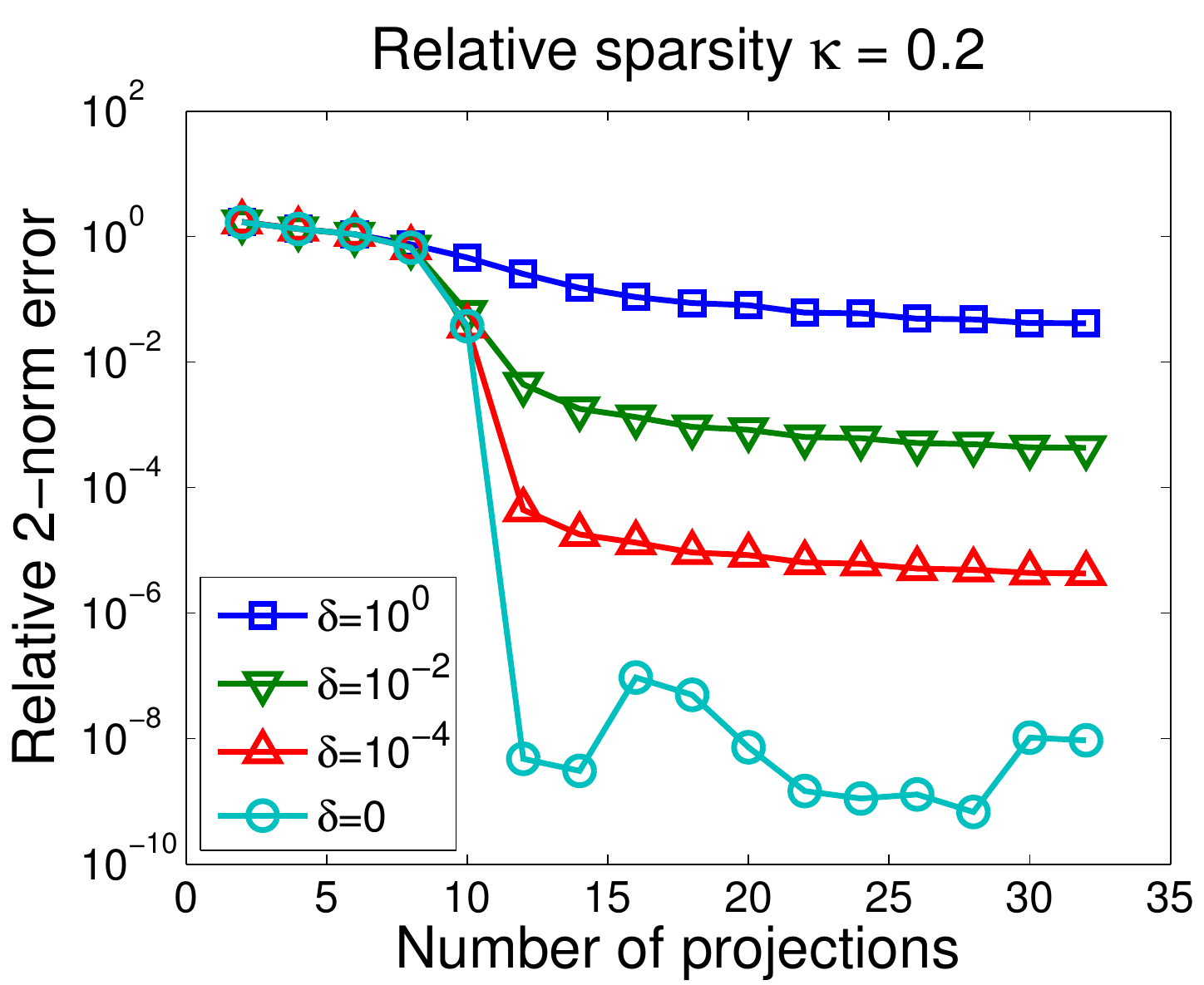}
 \includegraphics[width=0.345\linewidth]{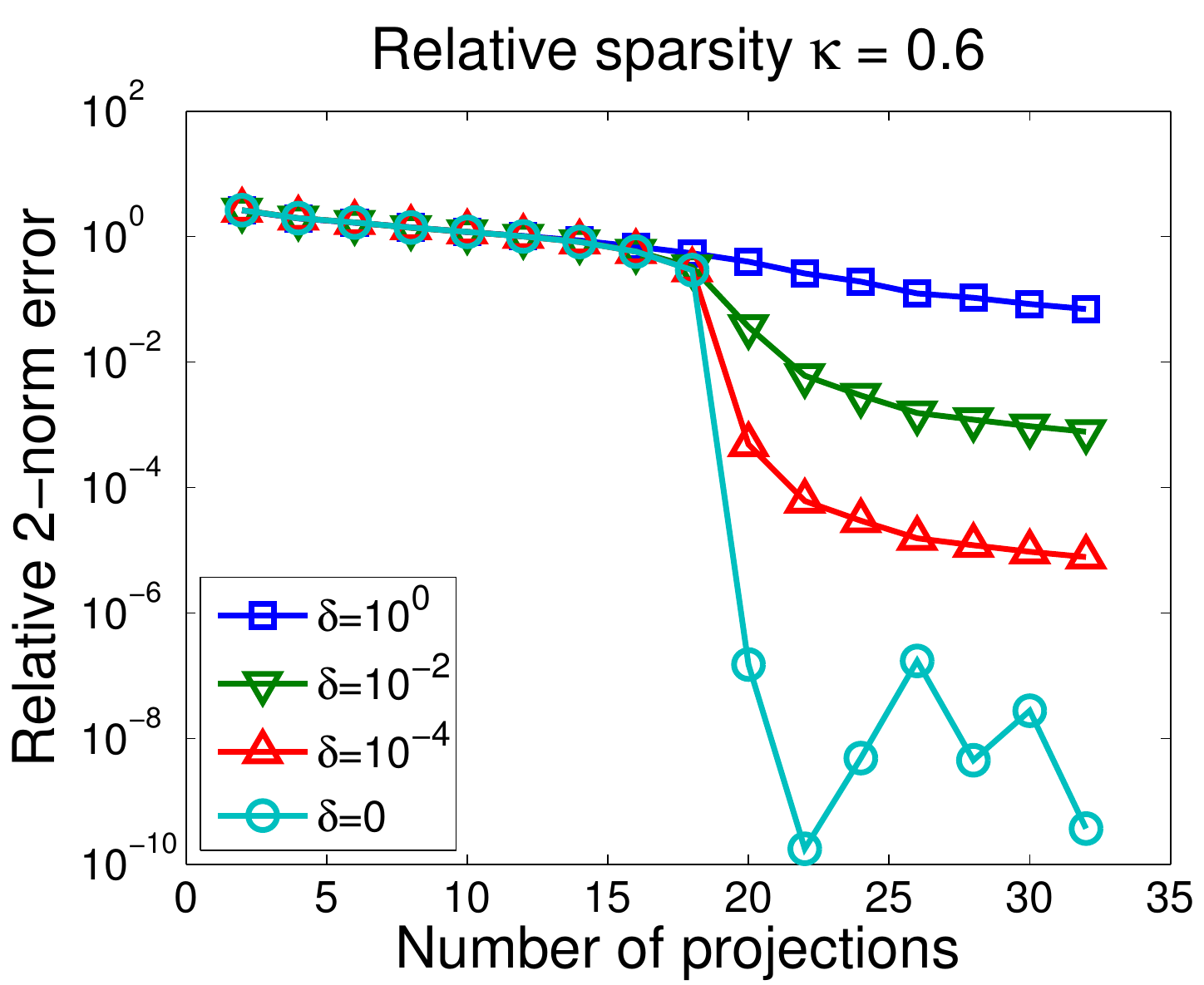}
 \includegraphics[width=0.29\linewidth]{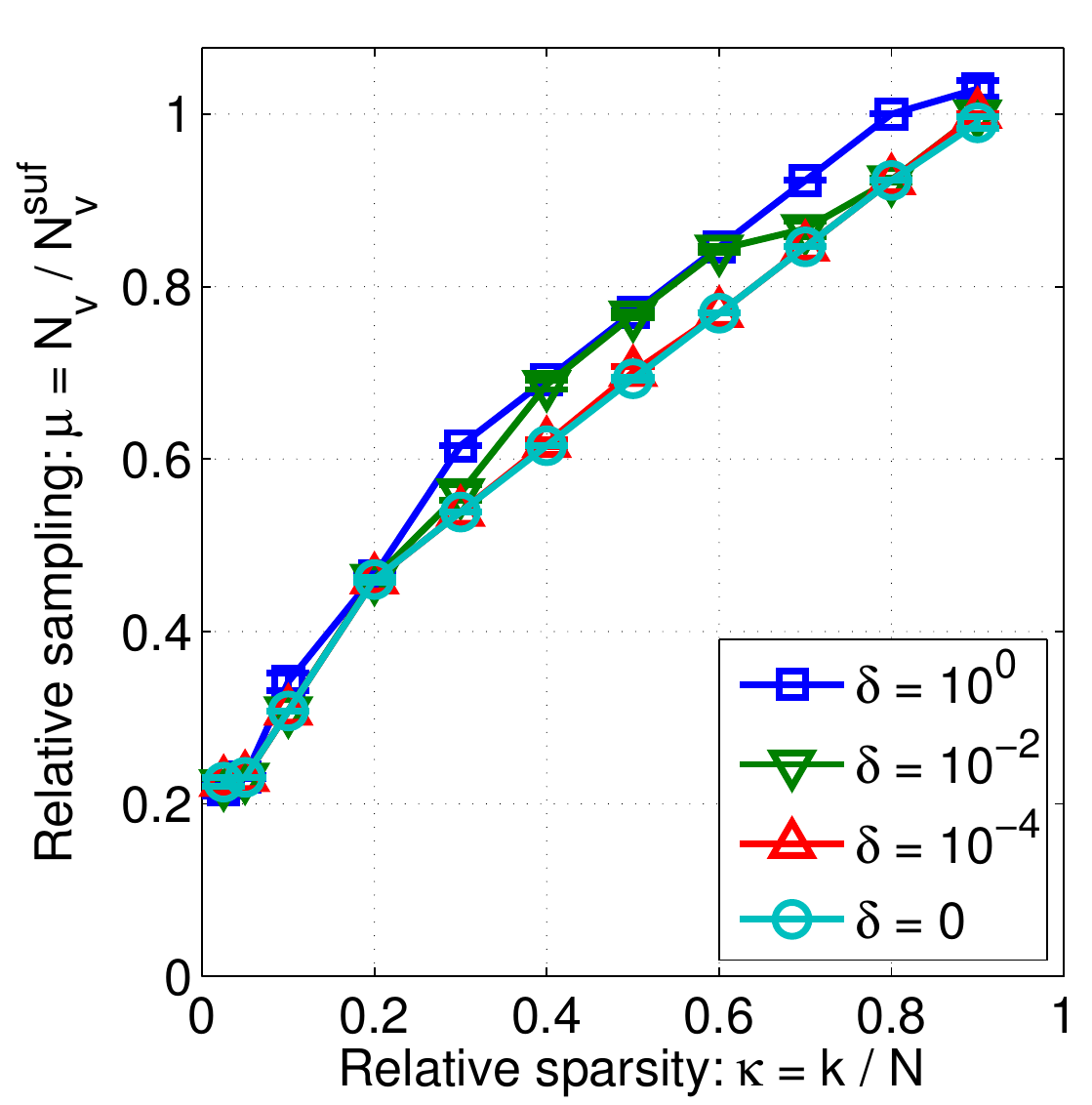}\\
\caption{Left, center: Relative errors of $\ellone_\ineqpar$ reconstructions vs.\ number of projections at different noise levels. Same \spikes{} image instances as in Fig.~\ref{fig:reconstruction_errors} with relative sparsity values $\relsparsity = 0.2$ (left) and $0.6$ (center). Right: Average sampling for $\ellone_\ineqpar$ recovery and $99\%$ confidence intervals for different noise levels. The ``$\delta=0$'' case is the noise-free \ellone{} result for reference.
\label{fig:noisy_reconstruction_errors}}
 \end{figure}

We consider the reconstruction problem
\begin{align} \label{eq:L1ineq}
\text{\elloneineq{}:}\qquad\qquad \min_\im  \ \normone{\im} \qquad
 \text{s.t.}\qquad \normtwo{\sysmat\im - \sino} \leq \ineqpar,
\end{align}
which can be solved in MOSEK by introducing a quadratic constraint:
\begin{equation}
 \argmin_{\im,w,t}\; \mathbf{1}^T w \quad \text{s.t.} \quad t = \sysmat \, \im - \sino \quad \text{and} \quad \|t\|_2^2 \leq \ineqpar^2 \quad \text{and} \quad -w \leq \im \leq w.
\end{equation}
We model each CT projection to have the same fixed x-ray exposure by letting the data in each projection  
$\sino_p$, $p=1,\dots,\nv$ be perturbed by an additive zero-mean Gaussian noise vector $e_p$ of constant magnitude $\normtwo{e_p} = \delta$, $p=1,\dots,\nv$. Hence, the noisy data are $\sino = \sysmat\imorig + e$, where $e$ is the concatenation of noise vectors for all projections. We use three noise levels, $\delta = 10^{-4}$, $10^{-2}$, $10^0$, corresponding to relative noise levels $\normtwo{e} / \normtwo{\sysmat \imorig}$ of $0.00016 \%$, $0.016 \%$ and $1.6 \%$. We reconstruct using \elloneineq{} with $\ineqpar = \normtwo{e} =  \sqrt{\nv} \cdot \delta$ and show the relative reconstruction errors from \eqref{eq:teststrongrecov} in Figure~\ref{fig:noisy_reconstruction_errors}. 

For $\delta = 10^{-4}$ and $10^{-2}$ the abrupt error drop when the image is recovered is observed at the same number of projections as in the noise-free case. The limiting reconstruction error is now governed by the choice of $\delta$ and not by the numerical accuracy of the algorithm as in the noise-free case. For the high noise level of $\delta=10^0$ no abrupt error drop can be observed. However, the reconstruction error does continue to decay after the number of projections for \ellone{} recovery seen at the lower noise levels and approach a limiting level consistent with 
the lower noise-level error curves.

In order to set up phase diagrams we must choose appropriate thresholds $\numthr$ to match the limiting reconstruction error at each noise level. In the noise-free case we used $\numthr=10^{-4}$ chosen to be roughly the midpoint between the initial and limiting errors of the order of $10^0$ and $10^{-8}$, respectively. Using the same strategy we obtain thresholds $10^{-2.5}$, $10^{-1.5}$, $10^{-0.5}$ for increasing noise level $\delta$. We determine the phase diagrams and for brevity and ease of comparison we only show the average sampling plots in Figure~\ref{fig:noisy_reconstruction_errors}. 
The low-noise phase transition is essentially unchanged from the noise-free case. With increasing noise level we see that the location of the transition is gradually shifted to higher $\relsampling$ values for the medium and large $\relsparsity$ values. At the high noise level and the largest $\relsparsity = 0.9$ we even see that there are instances that are not 
recovered (to the chosen threshold $\numthr$) at $\nvsuf$.

We conclude that the sparsity-sampling relation revealed by the phase diagram in the noise-free case is robust to low levels of Gaussian noise. For medium and high levels of noise, the phase diagram shows that a sharp transition continues to hold (for the particular noise considered) but the location of the transition changes to require more data for accurate reconstruction.

\section{Discussion} \label{sec:discussion}

Our simulation studies for x-ray CT show that for several image classes with different sparsity structure 
it is possible to observe a sharp transition from non-recovery to recovery in the sense that, on average, same-sparsity images require essentially the same sampling for recovery.
This is similar to what is observed in CS, but  
as explained in Section~\ref{sec:guarantees} no theory predicts that this should be the case for x-ray CT.
Based on the empirical evidence we conjecture that an underlying theoretical explanation exists.

The use of a robust optimization algorithm
limits the possible image size;
with MOSEK, we found $\imside > 128$ to be impractical.
Faster algorithms applicable to larger problems exist, however in our experience MOSEK is extremely robust in delivering an accurate \ellone{} solution, and the use of less robust optimization software may affect the decision of recoverability.
Furthermore, we found the relation between relative sparsity and relative sampling to hold independently of image size, so we expect that larger images can be studied indirectly through extrapolation.

The present work considers an idealized CT system, by focusing on recovery with \ellone{} rather than \elloneineq{}, and as such the quantitative conclusions are only valid for the specific CT geometry and image classes. Nevertheless we believe that our results can provide some preliminary guidance on sampling in a realistic CT system, as Section~\ref{sec:noise} indicates robustness of the relation between sparsity and sampling.

\subsection{Future work}

The phase diagram allows for generalization to increasingly realistic set-ups.
For example, more realistic image classes, sparsity in, e.g., gradient or wavelet domains, and other types of noise and inconsistencies 
can be considered by changing the optimization problem accordingly.
We expect that such studies will lead to improved understanding of the role of sparsity in CT.

Our earlier studies of TV reconstructions \cite{Joergensen_TMI:2013}
indicate a relation between the sparsity level in the image gradient and sufficient sampling for accurate reconstruction,
but due to the complexity of the test problems in that study we found
it difficult to establish any quantitative relation.
An investigation based on the phase diagram could provide more structured insight.
For instance, we might learn that TV reconstructions of a class of
``blocky'' or piecewise constant images exhibit a well-defined recovery-curve similar to
the ones in the present study.

With the present approach we always face the problem of possible non-unique solutions to \ellone{},
leading to phase diagrams that, in principle, depend on the particular
choice of optimization algorithm.
We expect that uniqueness of the \ellone{} solution can be studied by numerically verifying a set of necessary and sufficient conditions \cite{Grasmair2011}. We did not pursue that idea in the present work in order to focus on an empirical approach easily generalizable to other penalties, such as TV, for which similar unique conditions may not be available.

Finally, it would be interesting to study the in-class recovery variability,
i.e., why the \twopower{} class transition from non-recovery to recovery is more gradual.
Can differences be identified between instances that were recovered and ones that were not, e.g., in the spatial location of the non-zero pixels?
In \cite{PetraSchnoerr2014} it is found that the number of zero-measurements affects recoverability. In our case, the structure of \twopower{} leads to a higher and more variable number of zero-measurements, which might be connected with the larger in-class recovery variability.

\section{Conclusion}

We demonstrated empirically in extensive numerical studies a pronounced average-case relation between image sparsity and the number of CT projections sufficient for recovering an image through \ellone{} minimization. 
The relation allows for quantitatively predicting the number of projections
that, on average, suffices for \ellone{}-recovery of images
from a specific class, or conversely, to determine the maximal sparsity
that can be recovered for a certain number of projections.
  
The specific relation was found to depend on the image class with a smaller, but also more variable, number of projections sufficing for an image class of more structured images. Classes of less-structured images were found to exhibit a sharp phase transition from non-recovery to recovery.
We further demonstrated empirically that 
the sparsity-sampling relation
is independent of the image size and robust to small amounts of additive Gaussian noise.

With these initial results we have taken a step toward better quantitative
understanding of undersampling potential of sparsity-exploiting methods in x-ray CT.

\section*{Acknowledgments} 
This work was supported in part by Advanced Grant 291405 ``HD-Tomo'' from the European Research Council and by grant 274-07-0065 ``CSI:\ Computational Science in Imaging'' from the Danish
Research Council for Technology and Production Sciences.
JSJ acknowledges support from The Danish Ministry of Science, Innovation and Higher Education's Elite Research Scholarship.
This work was supported in part by NIH R01 grants CA158446, CA120540 and EB000225.
The contents of this article are
solely the responsibility of the authors and do not necessarily
represent the official views of the National Institutes of Health.

\bibliographystyle{AIMS}
\bibliography{sampling_eq_con_paper}

\end{document}